\documentclass[12pt]{amsart}
\usepackage{amsbsy}
\usepackage{graphicx,epsfig,subfigure,psfrag}
\begin{document}
%
%
\newtheorem{theorem}{Theorem}
\newtheorem{proposition}[theorem]{Proposition}
\newtheorem{lemma}[theorem]{Lemma}
\newtheorem{corollary}[theorem]{Corollary}
\newtheorem{definition}[theorem]{Definition}
\newtheorem{remark}[theorem]{Remark}
\numberwithin{equation}{section}
\numberwithin{theorem}{section}
\newcommand{\beq}{\begin{equation}}
\newcommand{\eeq}{\end{equation}}
\newcommand{\re}{{\mathbb R}}
\newcommand{\n}{\nabla}
\newcommand{\ren}{{\mathbb R}^N}
\newcommand{\iy}{\infty}
\newcommand{\pa}{\partial}
\newcommand{\fp}{\noindent}
\newcommand{\ms}{\medskip\vskip-.1cm}
\newcommand{\mpb}{\medskip}
\newcommand{\ssk}{\smallskip}
\newcommand{\BB}{{\bf B}}
\newcommand{\Am}{{\bf A}_{2m}}
\renewcommand{\a}{\alpha}
\renewcommand{\b}{\beta}
\newcommand{\g}{\gamma}
\newcommand{\G}{\Gamma}
\renewcommand{\d}{\delta}
\newcommand{\D}{\Delta}
\newcommand{\e}{\varepsilon}
\renewcommand{\l}{\lambda}
\renewcommand{\o}{\omega}
\renewcommand{\O}{\Omega}
\newcommand{\s}{\sigma}
\renewcommand{\t}{\tau}
\renewcommand{\th}{v }
\newcommand{\z}{\zeta}
\newcommand{\var}{\varphi}
\newcommand{\wx}{\widetilde x}
\newcommand{\wt}{\widetilde t}
\newcommand{\noi}{\noindent}
\def\com#1{\fbox{\parbox{6in}{\texttt{#1}}}}

\title
{\bf On 
 source-type solutions  and the Cauchy problem\\
for a doubly degenerate sixth-order thin film\\
 equation
   I. Local oscillatory
properties}

\author {M.~Chaves and V.A.~Galaktionov}

\address{Departamento de Matem\'atica Aplicada, Universidad de
Salamanca, Salamanca, SPAIN}
   \email{mchaves@usal.es}

\address{Department of Mathematical Sciences, University of Bath,
 Bath BA2 7AY, UK}
\email{vag@maths.bath.ac.uk}

\date{\today}




\begin{abstract}

As a key example, the sixth-order doubly degenerate parabolic
equation from thin film theory,
 $$
 u_t = (|u|^m |u_{xxxxx}|^n u_{xxxxx})_x \quad \mbox{in} \quad \re \times \re_+,
  $$
 with two parameters, $n \ge 0$ and $m \in (-n,n+2)$, is considered.   In this first part of the research,
 various local
 properties of its particular
 travelling wave and source type solutions are studied.
 Most complete analytic results on oscillatory structures of these solutions of changing sign
  are obtained for $m=1$ by an algebraic-geometric approach,
 with extension by continuity for $m \approx 1$.


\end{abstract}

\maketitle


\setcounter{equation}{0}
\section{Introduction: basic nonlinear model with complicated local
 and global properties of solutions}
\label{Sect1}
\setcounter{equation}{0}

\subsection{Higher-order degenerate parabolic PDEs: no potential,
monotone, order-preserving
 properties, not of divergence form, and no weak solutions}

We consider the {sixth-order
 parabolic  equation} from thin film theory,
 \beq
 \label{1}
 u_t = \big(|u|^m |u_{xxxxx}|^n u_{xxxxx}\big)_x \quad \mbox{in} \quad \re \times \re_+,
  \eeq
 with two parameters, $n \ge 0$ and $m \in (-n,n+2)$. This equation is {\em doubly
 degenerate} and
 contains the higher-order $p$-Laplacian nonlinearity
 $|u_{xxxxx}|^n$ (then  $p=2+n$),
 and another one
 $|u|^m$ of the porous medium type.  The
equation is written for solutions of changing sign, which is an
intrinsic feature of the Cauchy problem (the CP) with bounded
compactly  supported initial data $u_0(x)$ to be studied.

 The PDEs such as (\ref{1}), which are  called sixth-order
 {\em thin film equations} (the TFEs--6), were introduced by King
 \cite{Ki01} in 2001 among others
  for modelling of power-law fluids spreading on a horizontal substrate.
Equation (\ref{1}) is {\em quasilinear},  where the diffusion-like
operator includes  two parameters $m$ and $n$
 and {\em is not potential} (variational) and/or {\em monotone} in any functional
 setting and topology. It is not also an {\em operator of fully divergence form}, since the PDE
 admits just a single integration by parts, so that a standard definition of  weak
  solutions is entirely illusive.
  Of course, as a higher-order parabolic equation, (\ref{1})
   does not exhibit any {\em order-preserving} (via the Maximum Principle) features.

Overall, the TFE--6 (\ref{1}), as a typical example of a variety
of complicated nonlinear thin film models arising from modern
applications, represents a serious challenge to general PDE theory
of the twenty-first century, concerning principles and concepts of
understanding the common local and global features and properties
of its solutions, which also need proper definitions. For other
sixth-order TFEs including their derivation and mathematical
properties, see a survey in \cite{GBl6}, where further key
references are traced out.


 It is well-known that
{\em nonnegative solutions} of a wide class of higher-order TFEs
can be obtained by special
  non-analytic and often ``singular"
 $\e$-regularizations  and passing to the limit $\e \to 0^+$, that, in general,  lead to
free-boundary problems (FBPs). We refer to the pioneering work by
Bernis and Friedman \cite{BF1}
 and to the monograph on nonlinear parabolic PDEs \cite[Ch.~4]{Wu01}, where further references and
results can be found. Actually, such singular
$\e$-regularizations, as $\e \to 0^+$,  pose a kind of an
``obstacle FBP", where the solutions are obliged to be
non-negative by special free-boundary conditions that, in general,
are not easy to detect rigorously. Solutions of the CP cannot be
obtained by such techniques and require more involved
 and different analysis. In particular, {\em analytic}
 $\e$-regularizations, with $\e \to 0^+$, can be key for the CP,
  \cite{Gl4, GBl6}.

There is a large amount of pure, applied, and numerical
mathematical literature devoted to
 existence, uniqueness, and various local and asymptotic properties
of TFEs, especially, for the standard TFE--4:
 \beq
 \label{tfe4}
 u_t=- \big(|u|^n u_{xxx}\big)_x,
 \quad \mbox{where} \quad n>0.
  \eeq
 Necessary key references on various results of modern TFE theory
that are important for justifying principal regularity and other
assumptions on solutions will be presented below and, in
particular, are available in \cite{Gia08, Grun04} and  in a more
recent paper \cite{PetI}.
 See also \cite[Ch.~3]{GSVR}, where further references are given
and several evolution properties of TFEs (with absorption,
included) are discussed.
 However, even for simpler pure TFEs such as  (\ref{tfe4}),
 questions of local and global properties of solutions of the CP, their
 oscillatory, and asymptotic behaviour are not completely well understood
 or proved in
 view of growing complexity of mathematics corresponding to  {\em
 higher-order}
 degenerate parabolic flows.

Using this, rather complicated, and even exotic, model equation
(\ref{1}),
 we plan to explain typical and unavoidable difficulties
 that appear even in the study of local properties of compactly supported  solutions
 and their interfaces for higher-order degenerate nonlinear PDEs.
We then intent to give insight and develop some general
approaches, notions, and techniques,
 that are adequate and can be applied to a wide class of difficult degenerate parabolic (and
 not only parabolic) PDEs.
Since even the related ODEs of the fifth order for particular
solutions get very complicated with a higher-dimensional phase
space, we cannot rely on traditional ODE methods, which were very
successful in the twentieth  century for second-order ODEs,
occurred for many particular self-similar and other solutions,
with clear phase-planes.

 Instead, as a general idea, we propose to use
 parameter homotopy-continuity approaches using the fact that for
 some values of $m,\,n$ (e.g., for $m=1$ or $m=n=0$, etc.), the
 ODEs can be solved by some algebraic-geometric methods, or leads to easier linear equations.
  Then, we
 use a stable ``transversality-geometric" structure of the obtained solutions
to extend those into some surrounding  parameter ranges. However,
 global extension of those solutions are not straightforward at
 all and often we are obliged to apply careful numerical methods
 to trace out some solution properties and their actual existence.

  Therefore,  we are not restricted to {\em
 quasilinear}
 equations with semi-divergent operators. Without essential
 changes and hesitation, we may consider
 other
 {\em fully
 nonlinear}
 models such as the following formal
  parabolic PDE:
  \beq
  \label{FN1}
  |u_t|^{\s}u_t = \big(|u|^m |u_{xxxxx}|^n u_{xxxxx}\big)_x \quad \mbox{in} \quad
  \re \times
  \re_+ \quad (\s>0) \, ,
  \eeq
  where $\s=0$ leads to the quasilinear counterpart (\ref{1}).
  It is easy to propose other more artificial (and often awkward) versions of such PDEs without
  traces or remnants of  monotone, potential, divergence form, etc.,
  operators, which however can be used  for applying  our general mathematical concepts of
  analysis.


\subsection{On other models, results, and extensions}

 In the present first part,  in Sections \ref{Sect2}--\ref{SMap1},
  we present a  detailed  study of some, mainly local oscillatory and sign-changing,
   properties of
  {\em travelling wave} (TW) and {\em source-type solutions}.
   In this connection, let us mention the first important pioneering results on
  oscillatory source-type solutions of the fourth-order quasilinear parabolic equation
  of porous medium type (the PME--4) in the fully divergence form
  with the
   monotone operator in $H^{-2}$:
  \beq
  \label{m1}
  u_t =- (|u|^{m-1} u)_{xxxx}  \quad \mbox{in} \quad \re \times \re_+ \quad (m>1),
   \eeq
   which were obtained by Bernis \cite{Bern88} and in Bernis--McLeod \cite{BMc91}; see also
   \cite[\S~4.2]{Wu01} for further details, and \cite{GalRDE4n}
   for construction of a countable family of similarity solutions
   of (\ref{m1}).

 Notice that many key  features of a local oscillatory structure of solutions of
quasilinear degenerate PDEs with such operators do not essentially
change not only for equations of higher, $2m$th-order, but also
for similar odd-order {\em nonlinear dispersion equations} (NDEs).
In \cite{CGII}, as an illustration,
 we briefly
 review our approaches for the corresponding  fifth-order 
counterpart of (\ref{1}), which has the form (the NDE--5)
   \beq
 \label{1d}
 u_t = \big(|u|^m |u_{xxxx}|^n u_{xxxx}\big)_{x} \quad \mbox{in} \quad \re \times \re_+,
  \eeq
  Various
odd-order PDEs occur in nonlinear dispersion theory. As a key
feature, they exhibit finite interfaces, compacton behaviour, and
shock/rarefaction  waves; see first results in  Rosenau--Hyman
\cite{RosH93}, a survey in \cite[Ch.~4]{GSVR}, and \cite{GalNDE5},
as a more recent reference for shock wave behaviour for NDEs--5.
The counterpart of (\ref{1d}) for $m=n=0$ is the
 {\em linear dispersion} PDE
 \beq
 \label{5d}
 u_t = u_{xxxxx}\quad \mbox{in} \quad \re \times \re_+,
 \eeq
 which, by  natural continuity-homotopy issues,
  gives some clues about oscillatory properties of solutions for the nonlinear one (\ref{1d}).

It is curious that the performed local analysis of oscillatory
solutions close to finite interfaces remains similar for a number
of other PDEs. For instance, in this sense, the parabolic equation
(\ref{1}) has the counterpart which is a rather unusual PDE also
belonging to the class of  nonlinear dispersion equations
 \beq
 \label{d1}
 u_{tt} = \big(|u|^m |u_{xxxxx}|^n u_{xxxxx}\big)_{xx}.
  \eeq
 In its turn, a  ``local counterpart" of the nonlinear dispersion PDE
 (\ref{1d}) is the quasilinear sixth-order {\em hyperbolic}
 equation
  \beq
  \label{d2}
  u_{tt} = \big(|u|^m |u_{xxxx}|^n u_{xxxx}\big)_{xx}.
  \eeq

 The second part of our paper \cite{CGII} is devoted
  to a  global construction of source-type solutions and
  some general aspects concerning the Cauchy problem for
  (\ref{1}).
  As a first general rule therein,
 we show that, for correct understanding  the CP and some principal properties of
 its
  solutions, one should study   the limit $m,n
 \to 0^+$, where the  linear {\em tri-harmonic equation} occurs
  \beq
  \label{2}
  u_t= u_{xxxxxx}  \quad \mbox{in} \quad \re \times \re_+.
  \eeq
One can easily construct for (\ref{2}) the necessary TW and
source-type (ZKB) solutions, to be compared with solutions of its
nonlinear counterpart, the TFE (\ref{1}), at least for
sufficiently small $m$ and $n$. As a principal feature, we observe
that, for the CP, the solutions of both linear and nonlinear PDEs
are {\em oscillatory} near finite interfaces (for (\ref{2}),
  by implication,
  interfaces are at infinity, $x = \pm \iy$).

\smallskip

In higher-order quasilinear degenerate evolution PDEs of
parabolic, nonlinear dispersion, or hyperbolic types, the
identification of optimal local regularity properties of the FBP
and the Cauchy problem becomes very difficult, to say nothing
about general existence-uniqueness theory. Actually, we were
surprised to observe that many classic techniques of PDE theory
developed in the last fifty or so years not only cannot be applied
but also cannot provide us with necessary results {\em in
principle}, since these are strictly oriented to lower-order PDEs.

 Evidently,
related challenging difficulties appear even when constructing
standard particular self-similar solutions leading to higher-order
ODEs. These  cannot be studied by phase-plane analysis which is
usual for the first or second-order equations. Moreover, several
results, which have been obtained for  given third or fourth-order
PDEs by a complicated study of the topology of orbits, are often
hardly applied to similar ODEs of the order that is higher by one.
Many mathematical approaches are strictly attached to ODEs of the
given and sufficiently low order.

\smallskip

In this paper, we have specially chosen sixth- or fifth-order
models (\ref{1}) and (\ref{1d}) with complicated non-monotone and
non-potential operators to describe a certain general scheme to
study local and global properties of solutions including:

\smallskip

{\bf (I)} existence and uniqueness of travelling wave solutions, and

{\bf (II)} local oscillatory behaviour of solutions of changing sign near
interfaces.

\ssk

In a forthcoming paper \cite{CGII}, we continue this study of
(\ref{1}) and will concentrate upon:

\ssk

{\bf (III)} existence and uniqueness of source-type (ZKB-type) solutions,
and

\ssk

{\bf (IV)} general properties of solutions and proper setting of the
Cauchy problem.

\smallskip



\section{Travelling wave and source-type similarity solutions}
\label{Sect2}

Let us now introduce two classes of particular solutions to be studied in detail.

\subsection{Travelling waves}

These are simplest solutions of nonlinear PDEs of the form
 \beq
 \label{2.0}
u_{\rm TW}(x,t)= f(y), \quad y= x-\l t \quad (\l \in \re),
 \eeq
 where, on substitution into (\ref{1}), $f$ satisfies the autonomous ODE
 \beq
 \label{2.01}
 {\bf A}(f) \equiv |f|^m|f^{(5)}|^n f^{(5)}= -\l \, f.
 \eeq
 This is obtained on integration once with a zero constant via the required
 flux continuity at the interfaces, where $f=0$.

For $n=m=0$, the ODE (\ref{2.01}) is easy, with a linear bundle of
exponential solutions:
 \beq
 \label{2.02}
  f^{(5)}=-\l f \quad \Longrightarrow \quad
  f(y)={\mathrm e}^{\mu y}, \,\,\, \mbox{where} \,\,\,
  \mu^5=-\l.
   \eeq

\subsection{Source-type  solutions}
 \label{S2.2}

Equation (\ref{1}) admits the standard self-similar solutions
 \beq
 \label{2.1}
 u_{\rm S}(x,t)= t^{-\b} F(y), \quad y= x/t^\b,
 \quad \mbox{where} \quad
  \textstyle{
   \b = \frac 1{6+(m+n)},
    }
   \eeq
 so (\ref{2.1}) preserves the total mass of the solution:
 $$
 \mbox{$
 \int_\re u_{\rm S}(x,t) \,{\mathrm d}x = {\rm const.} \quad
 \mbox{for all} \,\,\, t \ge 0.
 $}
 $$
 Then $F$ solves a slightly different ODE again obtained on integration once with the zero constant:
  \beq
  \label{2.3}
{\bf A}(F) + \b F y \equiv  |F|^m  |F^{(5)}|^n F^{(5)}+ \b F y=0
\quad \mbox{in} \quad \re.
 \eeq
 In view of the symmetry group of scalings, if $F_1(y)$ is a solution
 of (\ref{2.3}), then
  \beq
  \label{2.4}
   \mbox{$
 F_a(y) = a^\g F_1\bigl(\frac y a\bigr), \quad \mbox{with} \quad \g= \frac 6{m+n},
  $}
 \eeq
 is a solutions for any $a>0$.
 Therefore, we are looking for a unit mass profile satisfying
  \beq
  \label{2.5T}
  \textstyle{
   \int F(y)\, {\mathrm d}y = 1.
   }
    \eeq
As another normalization, one can take the condition $F(0)=1$.

\smallskip

\noi{\bf Remark on finite propagation.} Finite propagation in the
ODEs such as (\ref{2.01}) and (\ref{2.3}), written in a semilinear
form (see below), i.e., existence of finite interface $y_0$ of
arbitrary sufficiently small solutions, was well-known for a long
time. We refer to the first ODE proofs in \cite{Bern88, BMc91},
\cite[p.~392]{Wu01}, and more involved energy estimates for
general related higher-order elliptic and parabolic  PDEs in
\cite{Bern01, Shi2} and survey in \cite{GS1S-V}. However, we must
admit that some extensions of energy methods to odd-order ODEs
such as (\ref{2.01}) or (\ref{2.3}), i.e., for nonlinear
dispersion operators involved, are not straightforward  and can be
technically rather involved.

 Thus, the principal question remained open is the oscillatory and
 non-oscillatory behaviour of solutions close to interfaces.

\smallskip

\subsection{Fundamental solution for $n=m=0$}

Then (\ref{2.1}) is the {\em fundamental solution}
 \beq
 \label{ffi1}
 b(x,t) = t^{-\frac 16} F_0(y), \quad y=x/t^{\frac16},
  \eeq
  of the
 tri-harmonic equation (\ref{2}). The corresponding  linear problem
 for the rescaled kernel denoted by $F_0$ reads
  \beq
  \label{2.51}
  \textstyle{
  F^{(5)}_0 + \frac 16 \, F_0 y=0 \quad \mbox{in} \quad \re, \quad \int
  F_0=1.}
   \eeq
  This has the unique solution $F_0(y)$ by classic linear theory,
   \cite{EidSys}; see \cite{CGII}  for extra properties of the kernel $F_0(y)$.
   Recall that we intend to use the linear rescaled kernel $F_0$ given by (\ref{2.51})
 in trying to understand the nonlinear one (\ref{2.3}), at least, for $m, \, n \approx 0$.

 \section{Local existence of positive solutions and maximal regularity}
\label{sectOsc}

 We now begin to study the behaviour of solutions of both ODEs (\ref{2.01}) and (\ref{2.3})
 near the interface point. For (\ref{2.01}), we just assume that the interface
 is at $y=0$ and set $f(y) \equiv 0$ for $y<0$. For
 (\ref{2.3}), assuming that
  $y_0<0$ is  the left-hand interface point of $F$, we perform the change $y-y_0
  \mapsto y$ and obtain equation (\ref{2.01}) with
    \beq
    \label{lll1}
   \l= \b y_0<0,
    \eeq
     up to  an exponentially small perturbation as $y \to y_0^-$.

  Thus, in both  cases, we  consider  (\ref{2.01})
 in a neighbourhood of the interface at
  $y=0$, where we scale out the parameter $|\l| \not = 0$, so that now $\l=\pm 1$:
  \beq
  \label{2.5}
   \mbox{$
   f^{(5)}= \mp |f|^{\a-1}f
  \quad \mbox{for} \quad y>0, \quad
 f(0)=0 \quad \big(\a= \frac {1-m}{1+n} \in (-1,1), \,\,\l=\pm 1 \big).
  $}
  \eeq
 Actually, the condition $\a \in (-1,1)$ holds in a wider parameter range
   \beq
   \label{2.5BB}
  n>-1, \quad m \in (-n,n+2).
   \eeq
The assumption $\a>-1$ ($m< n+2$) is purely technical
   that simplifies local analysis of $f(y)$ close to ``transversal" zeros
and makes it quite standard. For $\a \le -1$, some technicalities occur that
are often  not in the focus of the present study.

Recall that, in (\ref{2.5}), both $\l = \pm 1$ can occur  for the
TWs, and, always, $\l =-1$ for the source-type solutions.


We should clarify what kind of solutions we are looking for, and,
namely, which extra conditions are supposed to be posed at the
interface $y=0$ to create proper solutions of the CP. For
instance,
 one can look  for solutions
 $f \in C^4((-\d,\d))$, with a constant $\d >0$, so this demands
 \beq
 \label{2.5NN}
 f(0)=f'(0)=f''(0)=f'''(0)=f^{(4)}(0)=0,
  \eeq
but it is not clear whether  these correspond to the actual
regularity associated with the CP. Bearing in mind more
sophisticated ODEs that occur from models like (\ref{FN1}), we
would like to avoid any weak definitions of solutions  that are
based on integration by parts (for instance, for fully nonlinear operators, this
makes no sense).

As a hint, we can compare the desired regularity for the CP with
that for the standard zero height, zero-contact angle, zero
curvature, and zero-flux FBP, for which the free-boundary
conditions take the form
 \beq
 \label{FBPc}
f(0)=f'(0)=f''(0)=0 \quad \mbox{and} \quad (|f|^m|f^{(5)}|^n
f^{(5)})(0)=0.
 \eeq
Therefore,
   at the
interface $y=0$, solutions of the FBP satisfy for any $C>0$ and
some $C_1=C_1(C) \in \re$,
 \beq
 \label{FBr}
 f(y) = C y^3 + C_1 y^\g+... \,\,\, \mbox{for} \,\,\, y \ge 0 \,\,\,
 \Longrightarrow \,\,\,
 f \in C^2((-\d,\d)) \quad {\bf (\mbox{FBP})},
  \eeq
  which is true provided that
   \beq
   \label{FBP11}
   \mbox{$
   \g = \frac{8+5n-3m}{n+1}>3 \quad \Longrightarrow \quad
    m < \frac {5+2n}3.
    $}
  \eeq
Of course, (\ref{FBr}) exhibits  less regularity at $y=0$ than
(\ref{2.5NN}). In fact, the solutions of the FBP do not need
and/or  admit the zero extension for $y<0$
  (by the definition, this makes no sense for the FBP). If they do,
  this already corresponds to the CP.

For the CP, we need to use the concept of the {\em maximal
regularity}, which for the ODE (\ref{2.5}) simply means that setting $f(y)=0$ for $y<0$ yields
 \beq
 \label{mr}
 \mbox{solutions that are maximally smooth at $y=0$ admitted by the ODE}
  \quad {\bf (\mbox{CP)}}.
   \eeq
 The actual maximal regularity associated with the ODE under consideration
needs special local analysis close to $y=0^+$. It is worth
mentioning now that (\ref{2.5NN}) {\em are not} conditions of
maximal regularity in general, for arbitrary values of $m$ and
$n$.

\subsection{Existence, uniqueness, and nonexistence of
 positive TW solutions}

We begin by noting that, for $\l =-1$,
 equation (\ref{2.5}) admits the positive solution
  \beq
  \label{2.6}
  \textstyle{
  f_0(y)= \varphi_0 y^\mu, \quad \mbox{where} \quad \mu= \frac 5{1-\a}= \frac{5(n+1)}{m+n}
  \quad \mbox{and}
   }
  \eeq
  \vskip -.5cm
  \beq
  \label{2.61}
   \textstyle{
  \varphi_0= \bigl[ \frac{1}{\mu(\mu-1)(\mu-2)(\mu-3)(\mu-4))}
   \bigr]^{\frac{n+1}{m+n}}} \quad (\l =-1).
   \eeq
This formula makes sense in the intervals $\mu \in (4,+\infty)$,
$(2,3)$, and $(0,1)$, which give some relations between parameters
$m$ and $n$. For $\l=+1$, two intervals
 $\mu \in (1,2)$ and $(3,4)$ are accepted.
 Note that the regularity of the solution (\ref{2.6}) at $y=0$
 (with the trivial extension $f(y)\equiv 0$ for $y<0$) increases
 without bound as $\a \to 1^-$, since
  \beq
  \label{2.62}
   \mbox{$
   \mu= \frac 5{1-\a} \to + \infty \quad \mbox{as}
    \quad \a \to 1^-.
    $}
    \eeq
For instance,
 it follows that $f_0 \in C^4$ at the interface provided that
 \beq
 \label{mn1}
 \mbox{$
  m < \frac{n+5}4 \quad (\mu \in (4, +\infty)).
    $}
  \eeq
Below, we will explain the meaning of these positive solutions in
the Cauchy problem.

Notice that, for $\l =+1$, in the most convenient interval
(\ref{mn1}), the ODE (\ref{2.5}) does not admit solutions that are
strictly positive in an arbitrarily small neighbourhood of the
interface at $y=0$.

We first consider the natural ``smooth" version of solutions of the
maximal regularity, where the conditions (\ref{2.5NN}) hold, so
that $f_0(y)$ is such a solution in the interval (\ref{mn1}).
 The following classification of possible solutions holds.
 By {\em oscillatory solutions} of (\ref{2.5}) at the origin,
   we mean those that have infinitely many
   sign changes (for instance,
 {\em isolated} transversal zeros) in any arbitrarily small neighbourhood
$(0, \d)$ of the interface  at $y=0$.




\begin{proposition}
 \label{Pr.1}

 Let $(\ref{mn1})$ hold.  Then:

 \ssk

 \noi {\rm (i)} for $\l=+1$, any   $f(y) \not \equiv 0$ for $y \approx 0^+$ of the problem
  $(\ref{2.5})$,
   $(\ref{2.5NN})$
 is oscillatory at $y=0^+$, and

 \ssk

 \noi{\rm (ii)}  for $\l=-1$, there exists a unique local strictly monotone
 and positive solution of  $(\ref{2.5})$, which is given by
 $(\ref{2.6})$, and  all other solutions are  oscillatory as
 $y \to 0^+$.
  \end{proposition}

 \noi {\bf Remark: first comparison with the linear ODE.} The existence of a
1D manifold of positive solutions (with the interface parameter $y_0$) is
in good agreement  with the same conclusion for the linear equation
(\ref{2.02}), where the 1D manifold of positive exponentially decaying
functions  at the ``infinite interface" as $y \to -\infty$  is given by
 $$
 f(y)=
C\, {\mathrm e}^{ y}, \quad C>0.
 $$
   The rest of decaying
orbits are oscillatory as $y \to -\infty$. We thus prove that a
similar property persists for some $m, \, n>0$ in the interval
(\ref{mn1}).

\smallskip

  \noi{\em Proof of Proposition  $\ref{Pr.1}$.} (i) We argue by contradiction. Assume that $f(y)$ is a nontrivial solution
  of $(\ref{2.5})$ with $\lambda=+1$. Then, if $f(y)$ is positive in an interval $(0,y^*)$, it
follows from the equation (\ref{2.5})
  that $f^{(5)} < 0$ for $0<y<y^*$.
  After integration and taking into account the conditions (\ref{2.5NN}),
   we obtain that $f'(y) < 0$ and $f(y)<0$ on $(0,y_*)$, whence
   comes the contradiction. Similarly, there is no a negative
   solution on $(0,y_*)$.


\smallskip

(ii) The key idea of the proof  relies on the translation invariance
property of the equation. For $\lambda =-1$ and positive solutions,
consider equation (\ref{2.5}) in the form,
\begin{equation}
\label{2.5r}
 \mbox{$
 f^{(5)} = |f|^{\a-1}f \quad \mbox{for} \quad  y>0 \quad \big(\a= \frac{1-m}{1+n}<1\big).
 $}
\end{equation}
As a first consequence,  one can see after integrating this
equation  that any positive solution $f(y)$ of the ODE problem
(\ref{2.5}) is also convex and strictly increasing.

We next prove that, for $\lambda =+1$, the function $f_0= f_0(y)$
defined in (\ref{2.5})--(\ref{2.6}), is the unique positive
solution of the problem (\ref{2.5}), (\ref{2.5NN}). We begin by
showing that if $f_0(y)$ and $f(y)$ are different positive
solutions of (\ref{2.5}), then either they are ordered, or they
have mutually oscillatory behavior close to the origin (in the
sense that they intersect each other infinitely many times at
different values of $y$ approaching the origin, i.e., the
difference $f_0(y)=f(y)$ is oscillatory at $y=0$). Notice that, in
any domain of analyticity of $f(y)>0$ (where $f'(y) \not = 0$),
all intersections with the analytic $f_0(y)$ are isolated points.

We begin with the case $\alpha>0$. Assume that the assertion is
false and that there exists $y^*>0$ such that $f_0(y)$ and $f(y)$
are different functions ordered in $(0,y^*)$ and such that
$f_0(y^*)=f(y^*)$. If $f_0(y) < f(y)$ in this interval, it readily
follows from the equation in the new form (\ref{2.5r}) and the
regularity properties of the solutions that $f_0^{(5)}(y) \le
f^{(5)}(y)$ for $0<y \le y^*$. Hence, we obtain by integration in
$(0,y^*)$ that $f_0(y^*) < f(y^*)$, a contradiction with the
assumption. The same argument applies to the opposite inequality
and therefore, the assertion holds.

In order to establish the uniqueness result in both cases (ordered
solutions or with mutually oscillatory behavior), we apply similar
arguments to auxiliary solutions constructed by using the invariance
translation properties of the equation. Assume for contradiction that
there exists $y_0$ such that $f_0 (y)\neq f(y)$ with $f_0(y_0) < f(y_0)$.
Consider the auxiliary solution of (\ref{2.5}), (\ref{2.5NN}) defined as
follows:
 $$
 g(y)= \left\{ \begin{matrix} f(y - \e) \quad \mbox{for} \quad y
 \ge \e, \\ \quad\,\,\,\,
0  \quad \quad\,\,\, \mbox{for} \quad y
 < \e,
  \end{matrix}
  \right.
  $$
 where $\e>0$ is chosen small enough
such that $f_0(y_0) < g(y_0)$. So defined, we have that $g(y)$ is
a non-negative solution of $(\ref{2.5})$ and, clearly, it
satisfies by construction that $g(y) < f(y)$ for $y> \e$ close
enough to $\e$. Therefore we have by continuity and the previous
arguments that there exists a value $y^* \in (\e,y_0)$ such that
$g(y) \le f(y)$ in the interval $(0,y^*)$ and $g(y^*) = f(y^*)$.
Hence, it follows from equation $(\ref{2.5r})$ that $f_0^{(5)}(y)
\le g^{(5)}(y)$ for $0<y \le y^*$ and a contradiction at $y^*$
follows after integrating this inequality as above. A similar
contradiction argument applies if we assume that $f_0(y_0)>
f(y_0)$ by considering the auxiliary function $g(y) = f_0(y -
\e)$.

The uniqueness result for $\alpha <0$ is obtained by using similar
arguments but taking into account that in this case $h(s)=|s|^{\alpha-1}s$
is a decreasing function for $s>0$. This fact gives after subtracting and
integrating five times that
$$
 \mbox{$
f(y)-f_0(y) = \int\limits_0^y...\int\limits_0^s (f^\alpha-f_0^\alpha)(r) \, {\mathrm d}r\,...\,{\mathrm d}z,
 $}
$$ whence if the solutions are ordered close to the origin, the
signs of the left- and right-hand sides of this expression are
different, so the contradiction readily follows. The same argument
applies to $g(y)$ constructed as above, so mutually oscillatory
behavior is also disregarded in this case. This completes the
proof. $\qed$

\subsection{Positive solutions by fixed point theorem}
 \label{S3.2}

Indeed, we were lucky to have the {\em explicit} positive solution
(\ref{2.6}), which  is connected with the invariant scaling group
of transformations of the ODE (\ref{2.5}). We now sketch
 another approach
 to detecting a unique positive solution for more arbitrary
nonlinearities without using and relying on explicit calculus via
a scaling group.

We again consider the ODE (\ref{2.5r}) with $\l=-1$, i.e., look
for positive solutions of (\ref{2.5r}),
  where $f^\a$ can be replaced by more general functions $q(f)>0$
   for $f>0$. The main non-uniqueness difficulty for  (\ref{2.5r}) is
   that it always admits the trivial solution $f(y) \equiv 0$ if $\a>0$.
However, this disappears for
  the inverse function $y=y(f)$, for which
 $$
  \mbox{$
  f'= \frac 1{y'}, \quad f''= \big(\frac 1{y'}\big)'\frac 1{y'},
   \quad f'''=\big(\big(\frac 1{y'}\big)'\frac 1{y'}\big)'\frac 1{y'}, ... \, .
   $}
   $$
 Then we obtain the following ODE for $y(f)$:
  \beq
  \label{y1}
   \mbox{$
  f^{(5)} \equiv \big(\big(\big(\big(\big(\frac 1{y'}\big)'\frac 1{y'}
   \big)'\frac 1{y'}\big)' \frac 1{y'}\big)' \frac
  1{y'}\big)'\frac 1{y'} =f^\a \quad \mbox{for} \quad f>0.
   $}
   \eeq
Multiplying by  ${y'}$ and integrating over $(0,f)$ with
the zero boundary condition yields
 $$
  \mbox{$
\big(\big(\big(\big(\frac 1{y'}
 \big)'\frac 1{y'}\big)'\frac 1{y'}\big)' \frac
1{y'}\big)' \frac
  1{y'}= \int\limits_0^f f^\a y'_f \, {\mathrm d}f= f^\a y - \a \int\limits_0^f
 f^{\a-1} y \,{\mathrm d}f.
  $}
   $$
Multiplying again by ${y'}$ and integrating gives on the
right-hand side a smooth quadratic integral operator depending on
$y^2$, etc.

After five integrations like that, we obtain an integral equation
for $y(f)$ of the form:
 \beq
 \label{y4}
 y(f) = {\mathcal M}(y)(f)\quad \mbox{for} \quad f>0,
  \eeq
 with a smooth operator ${\mathcal M}$ of $y$ containing the maximum fifth degree $y^5$ polynomial
  dependence on $y$.
Therefore, this operator is a contraction in a space of continuous
functions provided that the typical integrals such as $\int_0
f^{\a-1}(\cdot)\, {\mathrm d} f$ converge, i.e.,
 \beq
 \label{y2}
  \a>0 \quad \Longrightarrow \quad m \in (0,1).
   \eeq
This restriction can be weakened by choosing special  classes of
functions $y(f)$ with prescribed envelopes as $y \to 0$.

Finally, Banach's Contraction Principle (see e.g.
\cite[p.~206]{KrasZ}) guarantees existence and uniqueness of a
positive solution of (\ref{y4}), and hence of (\ref{2.5r}). For
non-monotone changing sign solutions of (\ref{2.5r}), this inverse
function approach obviously fails. We will study such a behaviour
 in Section \ref{SS1}
by an extra scaling.

\subsection{Maximal regularity}

It is easy to see that (\ref{2.6}) describes the best (maximal)
regularity at the interface that is provided by the ODE
(\ref{2.01}). Indeed, (\ref{2.6}) established the only possible
balance between two terms in (\ref{2.01}). In other words, loosely
speaking, if such a balance is violated, the solution must behave
along the kernel (the null-manifold) of the nonlinear or linear
operator on both sides of (\ref{2.01}), which is obviously
trivial.
 Recall that this
maximal regularity reflects the intrinsic properties of the CP.
For the FBP, the regularity is different; cf. (\ref{FBr}).

\begin{proposition}
 \label{Pr.MR}
 Let $\alpha > 0$ and let the solution of $(\ref{2.01})$, $(\ref{2.5NN})$ be written as
  \beq
  \label{vv1}
   \mbox{$
   f(y)= y^\mu \varphi(y) \quad \mbox{for} \quad y \ge 0
   \quad \big( \mu= \frac 5{1-\a}= \frac{5(n+1)}{m+n}\big).
    $}
     \eeq
Then
$\varphi(y)$ is uniformly bounded for $y>0$.


 \end{proposition}

 \noi {\em Proof.}
 Let us  prove that $\varphi (y)$ is globally bounded, i.e., there
exists a positive constant $C$ such that
\begin{equation}
\label{upper}
 |\varphi (y)| \le C \quad \mbox{for all} \quad y >0.
\end{equation}
Assume that  $f\neq 0$, and, for every $y>0$, define
\begin{equation}
\label{sup}
 \textstyle{
 \bar y = \sup_{[0,y]}\,\,\{z \in [0,y]: \,\, |f(s)| \le
|f(z)| \,\, \forall s \in (0,y)\}.
 }
\end{equation}
Using the ODE such as (\ref{2.5r}) and integrating five times over
$(0,\bar y)$, we obtain by the definition in (\ref{sup}) that
$$
|f(\bar y)| \le C |f(\bar y)|^\alpha {\bar y}^5,
$$
so the desired inequality (\ref{upper}) readily follows for
$\bar y$ and, by (\ref{sup}), for every $y>0$. $\qed$

\smallskip

\subsection{Nonexistence of positive source-type similarity
profile}

 We finish this analysis by dealing with the question of
nonexistence of positive source-type solutions. We prove that, in
the basic interval (\ref{mn1}), the ODE problem (\ref{2.3}) has no
positive solutions with the symmetry conditions at the origin.

\begin{proposition}
\label{pr.2} There are no symmetric positive solutions with
compact support of the problem $(\ref{2.3})$, with conditions
$(\ref{2.5NN})$ for $F(y)$ at the interfaces.
\end{proposition}

\noi{\em Proof.} Let $F(y)$ be a positive solution of (\ref{2.3})
with compact support $[-y_0,y_0]$. It is clear from the equation
(\ref{2.3}) and the positivity assumption that $F^{(5)}(y)> 0$ in
$(-y_0,0)$.  By integrating this inequality over $(0,y_0)$
 and assuming the required regularity (\ref{2.5NN}) of $F(y)$ at $y=y_0$,
  we infer that  $F'(0) >0
$. Hence the symmetry condition fails that completes the proof. $ \qed$

\ssk

Thus, the sufficiently smooth source-type profile in the CP must
be oscillatory at the interfaces. We now begin to study the
character of such oscillations.


\section{Existence and uniqueness for the initial value problem}
 \label{sectIVP}

In order to analyze the existence and uniqueness problem, we
introduce the following (cf. \cite{Bern88}):

\begin{definition}
A function $f$ is said to be a weak $C^4$-solution of the equation $(\ref{2.01})$ on an
interval $I$, if $f\in C^4$ and for any $y_0,y \in I$,
\begin{equation}
\label{ivp1}
 \mbox{$
 f^{(4)}(y)-f^{(4)}(y_0)
= \pm \int\limits_{y_0}^y|f(s)|^{\alpha-1}f(s)\, {\mathrm d}s.
 $}
\end{equation}
\end{definition}

One can see that, for $\a\ge 0$, solutions are $C^5$ and hence
classical.
For $\a<0$, if we use the asymptotics  (\ref{vv1})  with uniformly
bounded (say periodic) oscillatory components  $\varphi$, weak
$C^4$-solutions exist for
 \beq
 \label{al14}
  \mbox{$
  \a > -\frac 14 \quad \big(\mbox{i.e.,} \,\,\,\frac{5\a}{1-\a} > -1\big).
   $}
   \eeq

 Consider the Cauchy problem for the equation (\ref{2.01}) with the
initial conditions
\begin{equation}
\label{ivp1c}
 f^{(j)}(y_0)= \alpha_j, \, \, j=0,1,2,3,4.
\end{equation}
It is clear that this initial value problem can also be written, by
introducing the function $h(y) = |f(y)|^{\alpha-1}f(y)$, in the equivalent
form,
\begin{equation}
\label{ivp}
\left\{%
\begin{array}{ll}
  f^{(4)}(y)= \alpha_4 \pm h(y), & \hbox{} \\
   h'(y)=|f(y)|^{\alpha-1}f(y), & \hbox{} \\
  f^{(j)}(y_0)= \alpha_j, \, \, \, j=0,1,2,3, & \hbox{} \\
  h(y_0)=0. & \hbox{} \\
\end{array}%
\right.    \\
\end{equation}
From the standard theory of ordinary differential equations, for
every $\alpha >0$, problem (\ref{ivp}) is solvable in some
neighborhood of $y=y_0$. Let $(f,h)$ be a solution. On the one
hand, it is clear that, for $\alpha \ge 1$, such solution is
unique since the Lipschitz condition on the nonlinearity  is
satisfied. On the other, for $\alpha < 1$, the only continuity is
guaranteed and in fact uniqueness does not hold, for instance, if
$\alpha_j=0$ for every $j=0,..,4$, as it is shown by means of the
construction of positive solutions in the previous section. A
partial answer to uniqueness is established in the following:

\begin{proposition}
Let $0< \alpha <1$. If in $(\ref{ivp})$
 $$
 \mbox{$
 \sum\limits_{j=0}^4 |\alpha_j| > 0,
  $}
  $$
  then the
initial value problem $(\ref{ivp})$ has at most one solution in a
neighbourhood of $y=y_0$.
\end{proposition}

\noi{\em Proof.} The proof follows the ideas in \cite{Bern88,
BMc91}. We prove uniqueness in a small right-hand  neighbourhood
of $y=y_0$. Uniqueness on the left is obtained in a similar way.
Denote by $j^*$  the smallest value of $j$ such that $\alpha_j$ is
non-trivial, and, for such $j=j^*$, introduce the function $$
  g_i(y)=
 \left\{
  \begin{matrix}\,\,
  \frac{f_i(y)}{(y-y_0)^j}, \, \, \,\quad \mbox{if} \,\,\, y>y_0, \\
  g_i(y_0)=0, \, \, \, \mbox{if} \,\, \,\,y=y_0.
  \end{matrix}
  \right.
$$
Assume for contradiction that there exist two different solutions $f_1$
and $f_2$ for small $y-y_0>0$. Due to the invariance
properties of the equation, we may assume without loose of generality that
$f_1 \ge f_2 \ge 0$. It is clear by integrating the equation and taking
into account the initial conditions that
$$
 \mbox{$
|f_1(y)-f_2(y)| \le \int\limits_{y_0}^y...\int\limits_{y_0}^s \big||f_1|^{\alpha-1}f_1(r) -
|f_2|^{\alpha-1}f_2(r) \big|\, {\mathrm d}r\,...\,{\mathrm d}z.
 $}
$$ Hence, after dividing by $(y-y_0)^j$, taking into account that
$h(s)$ is Lipschitz continuous away from $s=0$ and integrating
five times, we get $$ |g_1(y)-g_2(y)| \le
C(y-y_0)^{4+(1-\alpha)j}\,\max_{[y_0,y]}\, |g_1(s)-g_2(s)|. $$
Since $4+(1-\alpha)j >0$ for any $j \le 4$, the contradiction
follows for $y-y_0>0$ sufficiently small. The proof is complete.
$\qed$

\smallskip

Next we introduce a comparison principle which is valid for
$\alpha >0$ with $\lambda =-1$ (the positive sign on the
right-hand side in (\ref{2.5})) and $\alpha < 0$ with  $\lambda =
1$.

\begin{proposition}
Let $f_1$ and $f_2$ be two solutions of equation $(\ref{ivp1})$ on $[y_0,
\infty)$ satisfying $f_1^{(j)}(y_0) \ge f_2^{(j)}(y_0)$ for $j=0,1,2,3,4$, with
strict inequality for at least one of them. Then, $f_1 (y) > f_2(y)$ for
every $y \ge y_0$.
\end{proposition}

\noi{\em Proof.} The inequality is obvious by continuity and the
assumptions on the initial conditions on an small interval
$(y_0,y_0+\delta)$. We see that in fact the inequality holds for every $y
> y_0$. If the assertion is false, we denote by $y^* = \sup\{y \le y_0: \,
\, f_1(s) > f_2(s), \, \, s\in (y_0,y)\}$. It is clear that, at
this point, $f_1(y^*)=f_2(y^*)$. However, for $\alpha >0$ and
$\lambda =-1$, we obtain after subtracting and integrating five
times, $$
 \mbox{$
f_1(y^*)-f_2(y^*) = \sum\limits_{j=0}^4 (\alpha_1^j-\alpha_2^j) +
\int\limits_{y_0}^{y^*}...\int\limits_{y_0}^{s}
(|f_1|^{\alpha-1}f_1(r) - |f_2|^{\alpha-1}f_2(r))\,{\mathrm
d}r\,...\,{\mathrm d}y >0, $} $$ leading to a contradiction.  In a
similar way, a contradiction is obtained for $\alpha < 0$ and
$\lambda = 1$. $\qed$





\section{Oscillatory component: dissipative systems, periodic
behaviour, and heteroclinic bifurcations}
 \label{sectOsc2}

\subsection{Dissipative dynamical system for the oscillatory component}
 \label{SS1}

  We now study solutions of (\ref{2.5})
of changing sign. Bearing in mind (\ref{vv1}), we will be looking
for solutions in the form
 \beq
 \label{2.7}
  \mbox{$
 f(y) = y^\mu \varphi(s), \quad s = \ln y,
 \quad \mbox{where}
 \quad \mu = \frac 5{1-\a}= \frac{5(n+1)}{m+n}.
  $}
  \eeq
  Here $\varphi$ is called the {\em oscillatory component} of the
  solution, and is a sufficiently smooth function; see below.
The positive solution (\ref{2.6}) corresponds to the particular
case, where $\varphi(s) \equiv \varphi_0$ (see (\ref{2.61})), and
this exists in the corresponding parameter ranges. Besides
(\ref{2.6}), there exist many other changing sign solutions.

   According to
  Proposition \ref{Pr.MR},
     we
  know that, at least for $\a  \ge 0$, the oscillatory component satisfies
 \beq
 \label{2.7NN}
 \varphi(s) \quad \mbox{is bounded and $ \not \to 0$ as $s \to -
 \infty$.}
  \eeq
Therefore, (\ref{2.7}) gives a clear picture of typical regularity
of solutions at the interfaces,
 \beq
 \label{ii1}
  \begin{matrix}
  f \in C^4((-\d,\d)) \quad \mbox{for} \quad \quad \,\, m< \frac{n+5}4, \quad \,\, \smallskip\\
f \in C^3((-\d,\d)) \quad \mbox{for} \quad \frac {n+5}4 \le m<
\frac{2n+5}3,
\smallskip\\
f \in C^2((-\d,\d)) \quad \mbox{for} \quad \frac {2n+5}3 \le m<
\frac{3n+5}2.
\end{matrix}
 \eeq
 According to (\ref{FBP11}), in the first two ranges the CP always has a
 better regularity than the FBP.

Substituting the representation (\ref{2.7}) into the ODE
(\ref{2.5}), after simple manipulations via scaling properties,
yields the following autonomous
  equation for $\varphi(s)$:
   \beq
   \label{2.71}
 P_5(\varphi)=  \mp |\varphi|^{\a-1}\varphi
     \quad \mbox{in}
    \quad \re \quad (\l= \pm 1), 
     \eeq
     where $\a= \frac {1-m}{1+n}$.
 Here $P_5$ is an easily derived  linear differential operator of
 the form
   \beq
   \label{2.8}
 \begin{matrix}
     P_5(\var)=
{\mathrm e}^{(-\a\mu-1)s}({\mathrm e}^{-s}({\mathrm
e}^{-s}({\mathrm e}^{-s}({\mathrm e}^{-s}({\mathrm e}^{\mu
s}\varphi)')')')')'
  \smallskip\\
   \equiv    \varphi^{(5)}+5(\mu -2)\varphi^{(4)}
+5  (2\mu^2-8\mu+7)\varphi''' 
 \smallskip \\
 +\, 5 (\mu-2)(2\mu^2-8\mu  +5)\varphi'' 
 + (5\mu^4-40\mu^3+105 \mu^2 \smallskip \\
- \, 100\mu+24)\varphi' + \,
\mu(\mu-1)(\mu-2)(\mu-3)(\mu-4)\varphi.
 \end{matrix}
   \eeq

The best and simplest connection satisfying (\ref{2.7NN}) with the
interface at $s=-\infty$ is
 \beq
 \label{ps1}
\varphi= \varphi_*(s) \,\,\, \mbox{is a (non-constant) periodic
solution of (\ref{2.71})}.
  \eeq
On the other hand, any uniformly bounded solution of (\ref{2.71})
for $y \ll -1$ will fit.

We now concentrate on periodic connections.
 We list the
following properties that approach us to
 existence of a periodic
orbit of changing sign:

 \begin{proposition}
 \label{Pr.Per2}
 Let $\a \in [0,1)$ and $($cf. $(\ref{mn1}))$
   \beq
   \label{mn1NN}
    \mbox{$
   \mu= \frac 5{1-\a}= \frac{5(n+1)}{m+n} \ge 4, \quad \mbox{i.e.,}
   \quad m \le \frac{n+5}4.
    $}
    \eeq
 Then the fifth-order dynamical system $(\ref{2.71})$ satisfies:

 \ssk

\noi {\rm (i)} no orbits are attracted to infinity as  $s \to +
 \infty$, and

 \ssk

\noi {\rm (ii)} it is a dissipative system with  bounded absorbing
 sets defined  for any $\d>0$ by
 \beq
 \label{mn2}
  \mbox{$
  \limsup_{s \to +\infty} \, |\varphi(s)| \in
 B_\d^*=[-C_*-\d, C_*+\d], \quad \mbox{where} \quad  C_*=(5!)^{-\frac 1{1-\a}}.
 $}
  \eeq



\end{proposition}

The absorbing set in $\re^5$ includes all vectors
$(\varphi,\varphi',\varphi'', \varphi''',\varphi^{(4)})^T$, where
all derivatives are uniformly bounded by a certain
 constant.

\smallskip

\noi{\em Proof.} (i) The operator in (\ref{2.71}) is
asymptotically linear \cite[p.~77]{KrasZ} with the derivative at
the point  at infinity $P_5(\varphi)$ that has the characteristic
equation obtained by setting
  \beq
  \label{dd1}
  \varphi={\mathrm e}^{\rho s} \,\, \Longrightarrow \,\,
 p_5(\rho) \equiv (\mu+\rho)(\mu+\rho-1)(\mu+\rho-2)(\mu+\rho-3)(\mu+\rho-4)=0.
  \eeq
Therefore, all eigenvalues are real negative or non-positive:
 \beq
  \label{dd2}
   \mbox{$
 \rho_k=k-\mu \le 0 \quad \mbox{for}
 \,\,\, k=0,1,2,3,4, \quad \mbox{provided that $4-\mu \le 0$, or}
 \,\,\,
 m \le \frac {n+5}4.
 $}
  \eeq
Thus, $\varphi=\infty$ cannot attract orbits.

\smallskip

(ii) Actually, this follows from  (i) by using an extra scaling.
On the other hand, this is easy to see from the original equation
such as (\ref{2.5r}) for $\l=-1$ (for $\l=+1$ the proof is
precisely the same). Taking an arbitrary solution $f(y)$ defined
for all $y>0$, we integrate (\ref{2.5r}) five times to get that
 $$
  \mbox{$
  \bar f(y)= \sup_{z \in (0,y)} \, |f(z)|>0
   $}
  $$
 satisfies the inequality
  \beq
  \label{ii1N}
   \begin{matrix}
   \bar f(y) = C_0+C_1y+C_2 y^2+C_3 y^3 +C_4 y^4 + \int\limits_0^y\int\int \int \int
   |f(z)|^{\a-1} f(z)\, {\mathrm d}z
   \smallskip\smallskip\quad \qquad\\
    \le
    C_0+C_1y+C_2 y^2+C_3 y^3 +C_4 y^4   + \frac 1{5!} \, y^5 \, \bar
   f^\a(y),\qquad\quad
    \end{matrix}
    \eeq
    where $C_i>0$ are some constants.
 It follows that
  \beq
  \label{ii2}
   \mbox{$
  \bar f(y) \le C_0 + C_* y^{\frac 5{1-\a}} \quad \mbox{for all} \quad y \ge 0,
   $}
 \eeq
 provided that $\frac 5{1-\a} \ge 4$. In the variable (\ref{2.7}),
 this yields (\ref{mn2}).   $\qed$



\subsection{Periodic oscillatory component}

In connection with Proposition \ref{Pr.Per2}, existence of a
periodic orbit for dissipative dynamical systems is a standard
result of degree theory; see \cite[p.~235]{KrasZ}.
 Nevertheless, classic theory deals with non-autonomous periodic
 systems, with the given fixed period. Therefore, these
 results do not directly apply to the fifth-order dissipative autonomous
 dynamical system (\ref{2.71}), so we need some extra arguments to
 justify existence of periodic solutions for some values of $\a$ and how this disappears
when $\a$ gets negative.

It turns out that it is easier to establish first  where such periodic orbits are
nonexistent. The existence evidence will be presented later on.

\begin{proposition}
 \label{Pr.Non}
 {\rm (i)} The ODE $(\ref{2.71})$ with $\lambda=-1$ does not admit a nontrivial
$T_*$-periodic solution $\var_*$ for
 \beq
 \label{non.1}
  \mbox{$
  \mu \in (\frac 52, \mu_1^+) \,\,\, \mbox{or}
   \,\,\, \a \in (-1, -0.9655...),
   \quad \mbox{where}
  \quad
  \mu_1^- \approx 1.45608...\,, \,\,\, \mu_1^+ \approx 2.5439... \,.
   $}
 \eeq
{\rm (ii)} Analogously, the ODE $(\ref{2.71})$ with $\lambda=+1$ does not admit a
nontrivial $T_*$-periodic solution $\var_*$ for
 \beq
 \label{non.11}
  \mbox{$
 \frac 52< \mu \le \mu_*^+= 3.22474... \, , \,\, \mbox{or} \,\,
  \a \in (-1,\a_*^+), \,\,\, \mbox{where}
  \,\,\, \a_*^+ = \frac{\mu_*-5}{\mu_*} = -0.5505... \, .
   $}
 \eeq
 \end{proposition}

\begin{figure}
\centering
\includegraphics[scale=0.75]{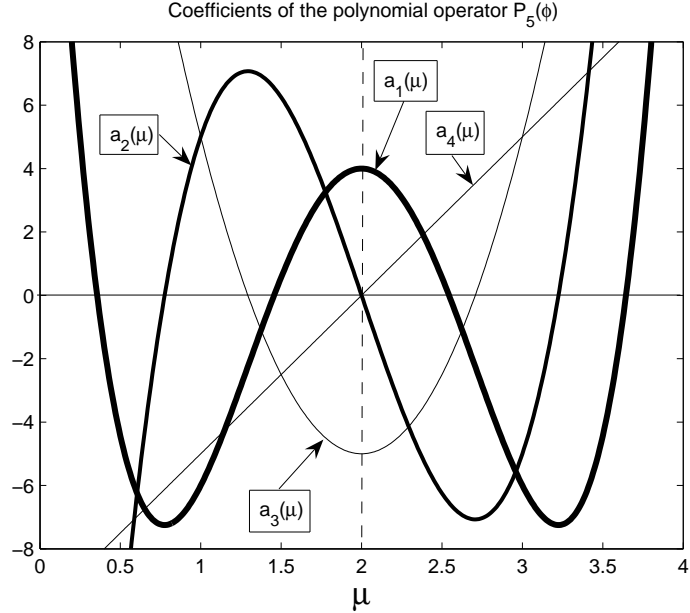}
\vskip -.3cm \caption{\small
 Coefficients of the polynomial operator $P_5(\var)$ in (\ref{2.8}) and
  (\ref{non.2}).}
   \vskip -.3cm
 \label{FCoeff}
\end{figure}

   \noi{\em Proof.} We write down (\ref{2.71}) as follows:
    \beq
    \label{non.2}
    P_5(\var)= \var^{(5)}+ a_4 \var^{(4)} + a_3 \var'''+a_2 \var''+a_1
    \var'+a_0 \var= \mp |\var|^{\a-1} \var \quad (\l= \pm 1),
     \eeq
     where the coefficients $a_i=a_i(\mu)$ are as in (\ref{2.8}).
     For convenience, we present graphs of these coefficients in
     Figure \ref{FCoeff}.

Multiplying (\ref{non.2}) by $\var$ and $\var'$
in $L^2(0,T_*)$
and integrating by parts yields,
 \begin{align}
 &
 \mbox{$
  a_4 \int(\var_*'')^2 - a_2 \int(\var_*')^2 + a_0 \int \var_*^2= \mp
 \int|\var_*|^{\a+1},
 $}
 \label{non.3} \\ \ssk
 &
  \mbox{$
  \int(\var_*''')^2 - a_3 \int(\var_*'')^2 + a_1 \int (\var_*')^2= 0,
  $}
 \label{non.4}
  \end{align}
    respectively. It is not difficult to see that those identities make sense for all $\a>-1$,
    so this restriction enters all our conditions.

    Note that multiplication by $\var''$  gives a new formal identity
    \beq
    \label{non.5}
    \mbox{$
    -a_4 \int(\var_*''')^2 + a_2 \int(\var_*'')^2 - a_0 \int (\var_*')^2=
  \pm
 \a \int|\var_*|^{\a-1}(\var_*')^2.
 $}
  \eeq
  However, on the right-hand side, there appears, after integration by parts,
 an integral with unknown converging properties, so using it could lead to wrong conclusions.
 More precisely, this integral makes sense for $\a>0$ only, which stays far away from the parameter ranges of interest.

Thus, here and later on various integrals of $\var_*$ and its derivatives
appear. The convergence of these integrals for $\a <0$ is checked by using
local properties of solutions of (\ref{2.71}). In particular, these are
easy for $\a>-1$, which justifies  identities (\ref{non.3}), (\ref{non.4}) in
the range of parameters under consideration, (\ref{2.5})--(\ref{2.5BB}).
Such convergence can be quite tricky for  negative
$\a<-1$ that are not under scrutiny in the present study.

\smallskip

 (i) To begin with our nonexistence purposes, it suffices to use just the
single second identity (\ref{non.4}). This gives that $\var_* \not ={\rm
const.}$ is nonexistent if
 \beq
 \label{non.6}
  a_1 \ge  0 \quad \mbox{and}
  \quad a_3 \le 0.
   \eeq
 It is easy to prove that there exist values  $\mu_1^-$ and $\mu_1^+$ roots of $a_1$ such that
$a_1 >0$ and $a_3<0$ on $(\mu_1^-,\mu_1^+)$, and hence, $\var_*$ does not
exist, at least, on such interval. Numerical approximation of these roots
shows that $\mu_1^- \approx 1.45608...$ and $\mu_1^+ \approx 2.5439...$\,.
The results are illustrated in Figure \ref{FCoeff}.

\smallskip

 (ii) In a similar way, the nonexistence result for $\lambda=1$ is extended to
(\ref{non.11}), by using identity (\ref{non.3}).
 This directly gives
that $\var_* \not ={\rm const.}$ does not exist if
 \beq
 \label{non.116}
  a_4 \ge  0  \quad \mbox{and}
  \quad a_2 \le 0,
 \eeq
which is true if $\mu \in ( \frac 52,\mu_*^+]$, where $\mu_*^+= 3.22474..$
is the largest root of $a_2$. $\qed$

 \smallskip

We continue  with establishing the principal fact that the
periodic solution $\varphi_*$  is  {\em hyperbolic} in some
parameter ranges of interest. This establishes an important
corollary concerning the type of {\em heteroclinic bifurcation},
and where it appears from; see below. This also reflects the
stability and instability properties of periodic solutions for
$\l=+1$ and $\l=-1$. We concentrate on the more important case
$\l=-1$, which in (\ref{2.71}) corresponds to oscillatory
behaviour of source-type solutions.

\begin{proposition}
 \label{Ptr.Stab}
 Let $\varphi_*(s)$ be a non-constant $T_*$-periodic solution of
 $(\ref{2.71})$, $\l=-1$, for
  \beq
  \label{ss.1}
   \mbox{$
   \mu \in (\mu_1^+,3], \quad \mbox{i.e.,} \quad \a \in (-0.9655...\,, -
   \frac 23].
    $}
    \eeq
Then $\var_*$ is hyperbolic.
 \end{proposition}



\noi{\em Proof.} For convenience, we write down (\ref{2.71}),
$\l=-1$ as \beq
  \label{AA.2}
 \var_*: \quad {\bf A}_-(\varphi) \equiv -P_5(\varphi)+|\varphi|^{\a-1}
  \varphi=0.
   \eeq
 Consider the eigenvalue problem for the
corresponding linearized operator
 \beq
  \label{AA.3}
 {\bf A}_-'(\varphi_*) \psi \equiv -P_5(\psi)- |\a| \, |\varphi_*|^{\a-1}
  \psi= \l_k \psi,
   \eeq
where we use that $\a<0$. Note that the potential here can be
rather singular at zeros of $\var_*(y)$ that can essentially
affect the setting of the eigenvalue problem and hence the actual
regularity of eigenfunctions. In particular, we always assume that
these zeros are transversal, so that the potential
$|\var_*(y)|^{\a-1}$ at zero, say, at $y=0$ is not singular as
$\sim \frac c{|y|^5}$ as $y \to 0$. Then by Hardy--Rellich-type
inequalities, we have got necessary embeddings that guarantee
compactness of the resolvent and hence discreteness of the
spectrum, which allows to deal with the hyperbolicity issue for
the periodic solutions.
   In further simple manipulations, we naturally assume that all these are justified.

Assume first, for simplicity, that an eigenvalue $\l_k$ is real.
Multiplying (\ref{AA.3}) by $ \psi$ and integrating over the
period $(0,T_*)$ by parts  yields
 \beq
  \label{AA.4}
  \begin{matrix}
-a_4 \int (\psi'')^2 + a_2  \int (\psi')^2 -a_0  \int \psi^2 -|\a|
\int |\var_*|^{\a-1} \psi^2
 =
\l_k \int
 \psi^2.
 \end{matrix}
  \eeq
  It follows that
   \beq
   \label{AA.5}
   \l_k <0 \quad ({\rm Im} \, \l_k =0), \quad \mbox{provided that}
   \quad a_4 \ge 0, \,\, a_2 \le 0, \,\, a_0 \ge 0.
   \eeq
   This gives the $\mu$-interval in (\ref{ss.1}).
   For $\l_k \in {\mathbb C}$, similarly,  first
   multiplying in $L^2$ by the complex conjugate $\bar \psi$, and
 next by $\psi$ the complex conjugate ODE,  summing up yields
  \beq
  \label{AA.6}
  \mbox{$
  -a_4 \int |\psi''|^2 + a_2  \int |\psi'|^2 -a_0  \int |\psi|^2
  -
|\a| \int |\var_*|^{\a-1} |\psi|^2
 =
   \frac {\l_k + \bar \l_k}2 \, \int
 |\psi|^2  < 0,
 $}
    \eeq
    i.e., ${\rm Re} \, \l_k <0$,
    under the same inequalities for the parameters as in
    (\ref{AA.5}). $\qed$

\ssk

The identities (\ref{AA.4}) and (\ref{AA.6}) can be used to detect
other hyperbolicity (and stability) parameter ranges of $\var_*$
and also for $\l=+1$.

As an important corollary, we obtain that, by classic
bifurcation-branching theory and implicit function theorem for
periodic solutions \cite[Ch.~6]{VainbergTr},  for $\l = - 1$, a
hyperbolic periodic solution $\varphi_*(s)$, existing at some
parameter value $\mu_0 \in (\mu_1^+,3)$, can be extended into a
small open neighbourhood of $(\mu_0-\d,\mu_0+\d)$. This implies
the following:

\begin{corollary}
 \label{Cor.1}
The $\mu$-parameter domain of existence of a periodic solution
$\varphi_*(s)$ of $(\ref{2.71})$, $\l=-1$ for $\mu \in (2,3)$, if
it is not empty, contains a connected interval $(\mu_{\rm h},
\mu_2)$, where $\mu_{\rm h} \ge \mu_1^+$ and $\mu_2 \le 3$.
 \end{corollary}




Finally, we expect that the periodic connection (\ref{ps1}) together with its
1D stable manifold as $s \to - \infty$
  is the only transition to the interface point at $s =-\infty$ ($y=0^+$),
 though this is difficult to justify completely rigorously.
 Also, we cannot prove  uniqueness (see below) of a periodic
solution of (\ref{2.71}).
  Before, the only result on existence of a periodic orbit
 is obtained in \cite[\S~7.2]{Gl4} for the third-order ODE like
 (\ref{2.71}) with $n=0$ and replacing the linear operator,
  $$
  P_5(\var) \mapsto P_3(\var)=\var'''+ 3(\mu-1) \var'' + (3 \mu^2- 6 \mu +2) \var'
 + \mu(\mu-1)(\mu-2) \var,
   $$
    by using the fact
 that this dynamical system is dissipative. Then classic theory
 \cite[\S~39.3]{KrasZ} applies together with a shooting-type
 argument. This results was improved in similar lines in
 \cite[\S~5]{PetI} with a sharper estimate on the periodic solution existence interval:
\beq
  \label{pp1}
   \mbox{$
  0<n<n_{\rm h} \in (\frac 32,n_{\rm +}), \quad \mbox{where}
   \quad n_{\rm +}= \frac 9{3+\sqrt 3}=1.9019238... \, .
    $}
    \eeq
    The actual heteroclinic bifurcation occurs at
 $$
    n_{\rm h}=1.7599... \quad (\mbox{the TFE--4}),
 $$
 and was calculated numerically, \cite{Gl4}.

\subsection{Heteroclinic bifurcations of periodic solutions}

 The
ODE (\ref{2.71}) is fifth-order quasilinear, is rather
non-standard, and has a non-Lipschitz nonlinearity. Moreover,
periodic solutions exist not for all $m,n>0$, namely, on some open
 (cf. Corollary \ref{Cor.1})
  interval $m \in (0,m_{\rm h}(n))$,
with the bifurcation exponent $m_{\rm h}$ to be discussed.
 The
classic results on existence of nontrivial periodic solutions
based on rotation vector field theory \cite[p.~50]{KrasZ} or
branching theory \cite[Ch.~6]{VainbergTr} do not apply to such
ODEs. Solutions that are bounded on the whole axis (which is fine
for (\ref{2.7NN})) also cannot be detected along the lines of
classic theory in \cite[p.~56]{KrasZ}. Moreover, we are interested
in solutions that are bounded (with non-zero limits) as $s \to
-\infty$ only. Several other techniques for existence of periodic
solutions also fail; see references and comments in
\cite[p.~140]{GSVR}.
 Therefore, we will rely
on careful  numerical evidence, especially when talking about the uniqueness of
periodic orbits and their stability.

\ssk

\noi{\underline{\sc The TW case $\l=+1$}.
In Figure
\ref{F1}, we show a stable periodic behaviour for the ODE
(\ref{2.71}) for $m=1$ and $n=0$ (so this is the standard thin film case) for $\lambda = +1$. It turns out
that such a periodic solution persists until a {\em
heteroclinic}-like bifurcation which occurs at
 \beq
 \label{mh0}
m_{\rm h}= 1.337968147... \quad (n=0, \,\, \lambda =+1).
 \eeq
 To compare with Proposition \ref{Pr.Non}, we present the corresponding universal critical
  values
  \beq
  \label{mh01}
   \mbox{$
\mu_{\rm h}= \frac{5}{1-\a_{\rm h}}=3.7370... \in (3,4), \quad \a_{\rm h}=
\frac{\mu_{\rm h}-5}{\mu_{\rm h}}=-0.3380... \, .
 $}
 \eeq
 In terms of the original parameters $m$ and arbitrary $n$, we have that
 the heteroclinic bifurcation occurs at
  \beq
  \label{mn771}
   \mbox{$
  m_{\rm h}(n) = \frac 5{\mu_{\rm h}}+ \big( \frac 5{\mu_{\rm
  h}}-1\big)n
  \approx 3.7370 + 2.3737 \, n.
   $}
  \eeq
Note that, for $\l=+1$, (\ref{mn771}) is essentially far from the
predicted nonexistence value (\ref{non.1}). It is interesting to
check whether the system (\ref{non.3}), (\ref{non.4}) (plus other
identities if any) can supply us in this case with a better
estimate of $\mu_{\rm h}$.

 This scenario of a heteroclinic
 bifurcation of stable and hyperbolic periodic solutions is
 typical for dynamical systems; see \cite[Ch.~4]{Perko}.
  Indeed, the hyperbolicity property
  as in Proposition \ref{Ptr.Stab}
   excludes  saddle-node bifurcations of periodic solutions (these demand existence of
  a  $\l \in {\mathrm i} \, \re$), at which $\var_*$ can disappear.
   Therefore, the appearance of a (stable) heteroclinic orbit as $\mu \to \mu_{\rm h}$
   in such dynamical systems is most plausible.

  Figure \ref{F2}  shows this typical
 formation of a
 heteroclinic orbit $-\varphi_0 \to \varphi_0$ as $m \to m_{\rm h}^-$.
 Here $\pm \varphi_0$ are constant equilibria of (\ref{2.71}) given
  as in (\ref{2.61}) by
   \beq
  \label{2.61NN}
   \textstyle{
  \varphi_0= \bigl[- \frac{1}{\mu(\mu-1)(\mu-2)(\mu-3)(\mu-4)}
   \bigr]^{\frac{n+1}{m+n}}} \quad (\l =+1).
   \eeq
Since by (\ref{mh01}),  for $n=0$ and $m=m_{\rm h}$, we have $\mu_{\rm h} \in (3,4)$,
both equilibria $\pm \varphi_0$ exist. These results have been
obtained by the {\tt MatLab} (the {\tt ode45} solver) with the
enhanced accuracy parameters Tols$\,=10^{-11}$ and the same
parameter of regularization in both degenerate and singular terms
in (\ref{2.71}).

 Notice that, in view of
(\ref{2.61}), solutions (\ref{2.7}) for $m \in( \frac 54, m_{\rm
h})$
 {\em do not} satisfy the last condition in (\ref{2.5NN}), i.e.,
 $$
 f^{(4)}(y) \,\,\, \mbox{is not bounded as $y \to 0^+$}.
 $$
 Nevertheless, the structure (\ref{2.7}) corresponds to the
 maximal regularity for the ODE (\ref{2.5}) and correctly establishes the balance of
  its nonlinear and linear terms.
In Figure \ref{F3}, we show a stable periodic behaviour for
$m=n=1$ and $\lambda=+1$.

\ssk

\noi \underline{\sc The source-type case $\l=-1$}. The case
$\lambda =-1$ that includes the source-type solutions, is shown in
Figure \ref{F4} for $m=n=1$. In this case, the periodic behaviour
is unstable as $s \to + \infty$ (as well as $s \to -\infty$),
though is clearly visible. For $n=0$, a ``heteroclinic"
bifurcation occurs at (see details in \cite[p.~142]{GSVR})
 \beq
 \label{mh1}
  m_{\rm h}= 1.909... \quad
(n=0, \,\, \lambda =-1).
 \eeq
 The corresponding other critical values are
  \beq
  \label{mh02}
   \mbox{$
\mu_{\rm h}= \frac{5}{1-\a_{\rm h}}=2.619... \in (2,3) \quad \mbox{and}
\quad \a_{\rm h}=
\frac{\mu_{\rm h}-5}{\mu_{\rm h}}=-0.909... \, .
 $}
 \eeq
Observe that this $\a_{\rm h}$ is sufficiently close to the
nonexistence one  $\a_*=-0.9655...$ in (\ref{non.1}) in
Proposition \ref{Pr.Non}. The corresponding critical values of
$m_{\rm h}(n)$
 are then given by
 \beq
  \label{mn772}
   \mbox{$
  m_{\rm h}(n) = \frac 5{\mu_{\rm h}}+ \big( \frac 5{\mu_{\rm
  h}}-1\big)n
  \approx 1.909 + 0.909 \, n.
   $}
  \eeq

\begin{figure}
\centering
\includegraphics[scale=0.65]{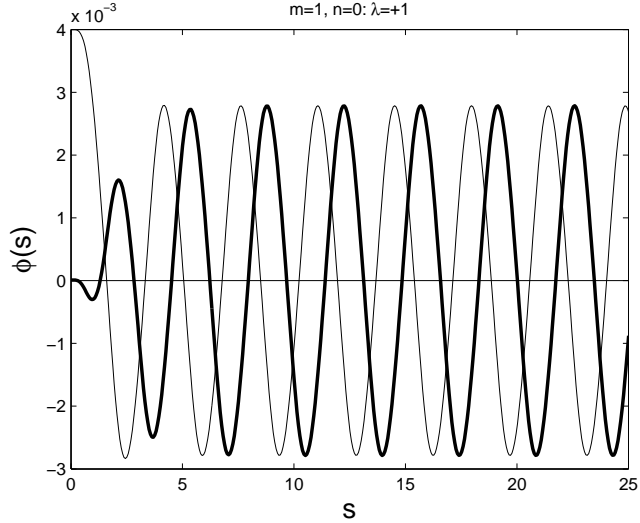}
\vskip -.3cm \caption{\small
 Stable periodic behaviour for the ODE  (\ref{2.71}) for $m=1$, $n=0$, $\lambda =
 +1$.}
   \vskip -.3cm
 \label{F1}
\end{figure}

\begin{figure}
\centering
\includegraphics[scale=0.65]{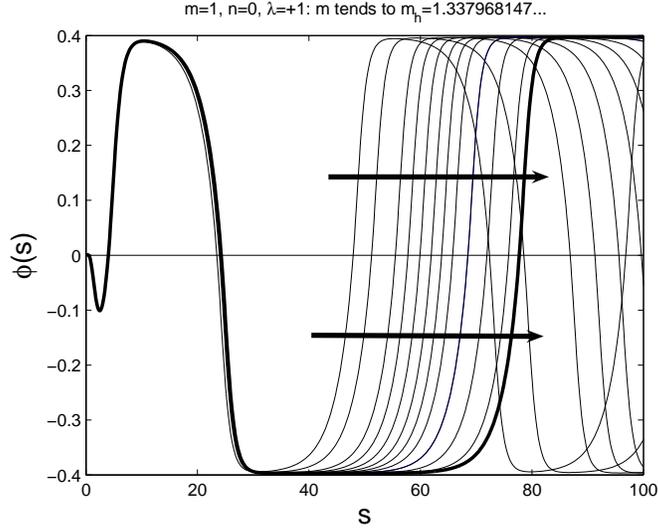}
\vskip -.3cm \caption{\small Formation of a heteroclinic orbit as
$m \to m_{\rm h}^-=1.3380...$ for $n=0$, $\l=+1$.}
   \vskip -.3cm
 \label{F2}
\end{figure}

\begin{figure}
\centering
\includegraphics[scale=0.65]{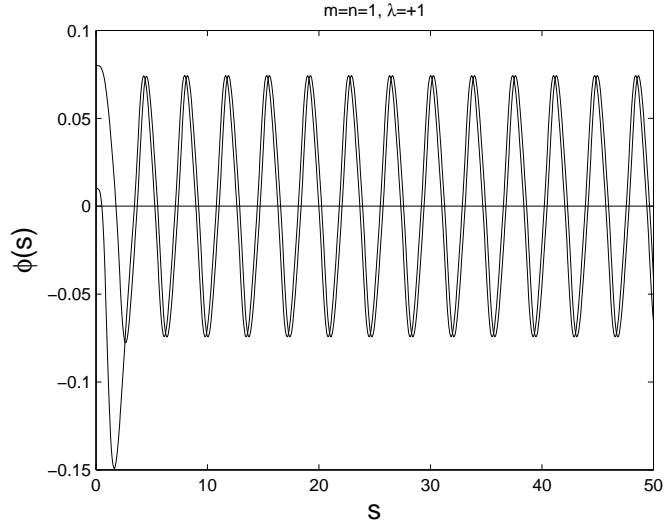}
\vskip -.3cm \caption{\small
 Stable periodic behaviour of (\ref{2.71}) for $m=n=1$, $\lambda =
 +1$.}
   \vskip -.3cm
 \label{F3}
\end{figure}

\begin{figure}
\centering
\includegraphics[scale=0.65]{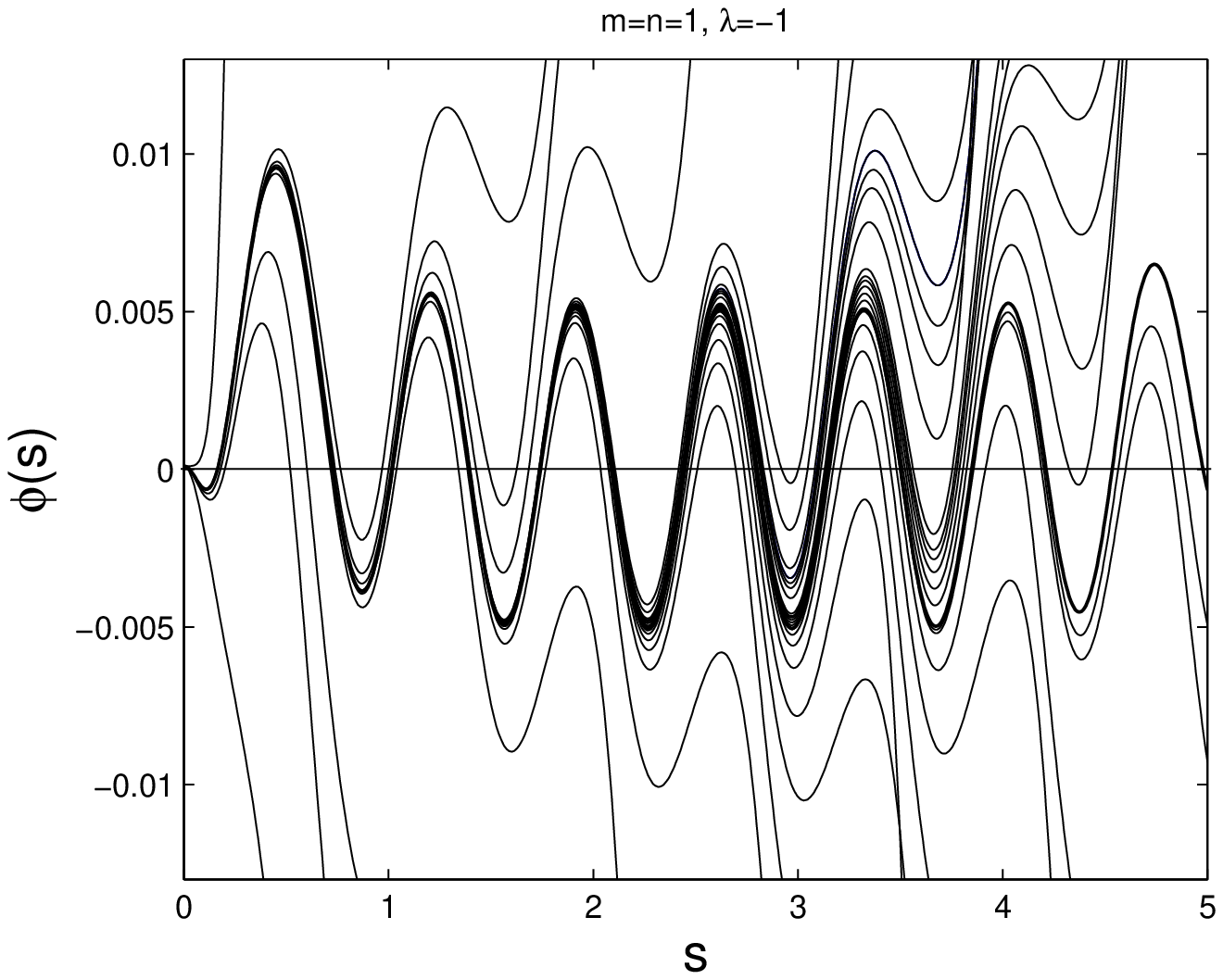}
\vskip -.3cm \caption{\small
 The trace of an unstable periodic orbit of (\ref{2.71}) for $m=n=1$, $\lambda =
 -1$.}
   \vskip -.3cm
 \label{F4}
\end{figure}

Let us
 state the following conjecture on periodic solutions and
their stable manifolds.

\smallskip

\noi {\bf Conjecture \ref{sectOsc2}.2.} {\em For all $n \ge 0$,
there exists $m_h(n)$ given by $(\ref{mn771})$  for $\l=+1$
$($resp. by $(\ref{mn772})$ for $\l=-1)$, such that for all $m \in
(0,m_{\rm h}(n))$,  the ODE $(\ref{2.71})$ with $\l = +1$ $($resp.
$\l=-1)$:

\ssk

 \noi{\rm (i)} has a unique stable $($unstable$)$ as $s \to +\infty$
periodic solution $\varphi_*(s)$ of changing sign;

\smallskip

\noi{\rm (ii)} as $s \to -\infty$, the  periodic solution
$\varphi_*(s)$ is unstable and  has  a 1D stable manifold; and

\ssk

 \noi{\rm (iii)} there exists a heteroclinic bifurcation at  $m_{\rm
 h}(n)$, so $\var_*(s)$ is nonexistent for $m \ge
 m_{\rm h}(n)$}.


\smallskip

 It is worth mentioning the obvious instability of
$\varphi_*(s)$ as $s \to -\infty$ (while approaching the interface
point at some $y=y_0$) is associated with the translational
invariance of the ODEs such as (\ref{2.5}) admitting shifting in
$y_0$. The precise meaning and the significance of the conclusion
(ii) will be explained in the second part of the paper
\cite{CGII}, which will be key for a well-posed shooting of a
global source-type similarity profile.





\subsection{Oscillatory component: nonlinear dispersion  model}

Substituting TW solutions (\ref{2.0}) into (\ref{d1}) and
integrating twice with zero constants of integration yields the
same ODE (\ref{2.01}) with the only change on the right-hand side:
 $$
 -\l \mapsto \l^2>0.
  $$
  Therefore, for the oscillatory component we obtain the same
  equation (\ref{2.71}), and in examples and figures presented
  there, we always take $\l=-1$, so that $-\l=1>0$. In particular, the
  oscillatory behaviour near interfaces persists until the
  bifurcation exponent (\ref{mh1}) for $n=0$.


\section{Existence and uniqueness of oscillatory TW solutions for $m=1$}
 \label{S6}

In this section, we analyze the problem of existence and
uniqueness of the oscillatory travelling waves solutions, as well
as the character of their oscillations. We begin by proving
existence of the oscillatory TW solutions via a different
approach.

\subsection{Existence and uniqueness: an algebraic-geometric approach to periodic orbits}
 \label{S6.1}

Here, we develop an alternative approach to existence of periodic orbits
of (\ref{2.71}) for both $\lambda$, positive and negative. We focus our
attention on the construction of solutions with right-hand side
interfaces, $f^+(y)$ and $f^-(y)$ corresponding to $\lambda
>0$ and $\lambda < 0$ respectively. By
reflection, we obtain solutions with left-hand
side interfaces, namely $f^+(-y)$ for $\lambda <0$ and $f^-(-y)$ for
$\lambda >0$.

We begin with  the original equation for $f(y)$, (\ref{2.5r}),
which, for the case  $m=1$ (and putting for convenience  $|\l|$ to
make calculations easier), takes especially simple  ``weakly
nonlinear" form
 \beq
 \label{b1}
 f^{(5)}= 5!\, {\rm sign} \, f.
  \eeq
  Let us explain the main ingredients of the strategy of our
  construction.

By using the properties of the solutions, we may consider for
convenience, that the positivity domain of $f_0$ is given by
$(-1,0)$, so $f(0)=f(-1)=0$ and $f(y)>0$ on $(-1,0)$. Hence, such
a function denoted now by $f_0(y)$ satisfies on $(-1,0)$,
    \beq
 \label{b2}
  \mbox{$
 f^{(5)}=5! \quad \Longrightarrow \quad f_0(y)= y(y+1)[a + by+ cy^2+
 y^3],
  $}
   \eeq
where the constants $a,b,c$ are chosen so that $f_0(y)>0$ on
 $(-1,0)$.

We next extend $f_0(y)$ to $y>0$ as follows in three steps:

\ssk

\noindent (i) take $-f_0(y)$,

\ssk

\noindent (ii) shift it to the right  by $y_1=1$ to get $-f_0(y-1)$, and

\ssk

\noindent (iii) rescale it by  the invariant scaling group
   for (\ref{b1}) to get for $y>0$,
    \beq
    \label{b3}
     \mbox{$
    f_0(y) \equiv  -G^5 f_0 \bigl(\frac y{G} - 1\bigr), \quad \mbox{with a parameter $G>0$}.
     $}
     \eeq

In order to have a smooth solutions on $[-\d,\d]$, we need four matching
conditions on zero jumps of  four consecutive derivatives,
 \beq
 \label{b4}
[f_0'](0)=[f_0''](0)=[f_0'''](0)=[f_0^{(4)}](0)=0.
 \eeq
Thus, we obtain {\em four algebraic equations} with four parameters
$a,\,b,\, c$, and $G$. They can be written as:
\beq \label{S}
\left\{%
\begin{array}{ll}
  (1-G^4)a + G^4b -G^4c =-G^4, & \hbox{} \\
   (1+G^3)a + (1-2G^3)b+ 3G^3c=4G^3, & \hbox{} \\
  (1+G^2)b+(1-3G^2)c=-6G^2, & \hbox{} \\
 (1+G)c=4G-1. & \hbox{} \\
\end{array}%
\right.
\eeq

We remark that by the argument of continuation of $f_0$  to a solution of
the problem with a finite right-hand  interface, the condition
 $$
 0<G<1
 $$
  is required.
Therefore, we restrict ourselves to analyzing system (\ref{S}) in this
range of parameters.

\ssk

We begin by noting that, for every fixed value of $G \in (0,1)$,
(\ref{S}) is a linear system in $a,b,c$ of four equations.
Considering the three last equations, one can easily  check  that
the determinant of the associated matrix
 $$
 d=(1+G^3)(1+G^2)(1+G)
   $$
    is not zero, and, hence, for every fixed $G \in (0,1)$,
the system of the last three equations has a unique solution. Therefore,
it will be the solution of the complete system (\ref{S}) if the
determinant of the matrix,
\begin{equation}\label{M}
\begin{pmatrix}
1-G^4 &  G^4 & -G^4 & -G^4  \\
1+G^3 & 1-2G^3 & 3G^3 & 4G^3 \\
0 & 1+G^2 & 1-3G^2 & -6G^2  \\
0 & 0 &  1+G & 4G-1
\end{pmatrix}
\end{equation}
is zero. Denoting by $D(G)$ the determinant of this matrix, we
have that $D(0)$ and $D(1)$ are both negative. Hence, it readily
follows that either no roots of $D(G)$ are available, or at least
two roots exist. After analyzing in more detail $D(G)$, one can
see that there exist two roots, that can be obtained numerically:
 $$
 G_1 =  0.178318... \quad \mbox{and} \quad
 G_2 = 0.7060378...\,.
  $$
However, the positivity condition of $f_0$ in $(-1,0)$ implies
$a<0$ or, equivalently, by using the first equation in (\ref{S}),
gives
 $1+b-c>0$. Manipulating the two last equations in (\ref{S}),
one obtains that this restriction holds, if only if,
 $$
 F(G)=3G^2-10G+3>0.
   $$
    One
can check that this is satisfied by $G_1$ and allows to  disregard $G_2$.

It is key that the existence of the value $G_2$ also plays an
important role in the analysis, since it provides, by means of the
construction explained above, solutions of the equation with the
left-hand interface for $\lambda <0$, i.e., \beq
 \label{b2NN}
 f^{(5)}= -5!\, {\rm sign} \, f.
\eeq

These arguments allow to state the following:

\begin{theorem}
 \label{Th.Exm1}

 {\rm (i)} The above algebraic system $(\ref{S})$ has a unique solution $G_1 \in (0,1)$
  satisfying the
 positivity condition of $f_0=f_0(y;G_1)$ in $(-1,0)$ and a unique solution $G_2\in (0,1)$
  satisfying the negativity
 condition of $f_0=f_0(y;G_2)$ in $(-1,0)$.

 \ssk

 {\rm (ii)} Equation $(\ref{2.71})$, with the  sign $``+"$ $($respectively $``-")$
  and $\a=0$, has a nontrivial periodic
 solution $f^+(y)$ $($respectively $f^-(y))$, with the period
 \beq
 \label{per1}
     T_*= 2 \ln G_1 \quad \big( \mbox{resp.} \quad T_*=2 \ln G_2 \big).
     \eeq
 \end{theorem}

 \noi{\em Proof.} (i) We prove the result for $\lambda$ positive
 and the solution corresponding to $G_1$. For $\lambda
 <0$ the same arguments apply.

Once a solution of the algebraic system  (\ref{S}) has been found,
we can extend the solution $f_0(y)$ similarly to the next interval
of oscillations, etc. indefinitely. Then the sequence of positive
and negative humps will converge according to some geometric
series. For instance, the right-hand interface of the interval of
every extension, given by the standard geometric series, $$ b_n= G
+ G^2 + G^3 + ...+ G^n \quad (G=G_1 \,\,\, \mbox{or} \,\,\, G_2),
$$ converges as $n \to \iy$ to the right-hand interface of the
solution, $\frac G {G-1}$.

It is not hard to see that, after reflecting and moving the interface to
$y=0^+$, this gives precisely the
behavior (\ref{2.7}), i.e., 
 \beq
 \label{b5}
  f_0(y) = (y_0-y)^5 \varphi_*(s), \quad s= \ln (y_0-y),
   \eeq
 where $y_0 = \frac G{G-1}$ denotes the right-hand interface of $f_0$.
 We next prove that the oscillatory component $\varphi_*$ is periodic with
   the period (\ref{per1}). In fact, we show that $\varphi_*$
 satisfies a stronger property,
  $$
  \var^*(s+\ln G)= -\var^*(s)\quad \mbox{for every $s \in
 \mathbb{R}$}.
  $$

 Hence,
  $$
  \var^*(s+2\ln G)= -\var^*(s + \ln G)= \var^*(s)
   \quad \mbox{for every $s\in
\mathbb{R}$},
 $$
  and the periodicity stated for $\var^*$ follows.
According to the definitions in (\ref{b5}) and  the expression of
the right-hand interface $y_0$, one obtains that
 $$
   \var^*(s+\ln G)= G^{-5}{\mathrm e}^{-5s}f_0(y_0-G{\mathrm e}^s) \quad \mbox{and} \quad
   \var^*(s)={\mathrm e}^{-5s}f_0(y_0-{\mathrm e}^s).
 $$
 Assume, for instance, that $s$ is such that $y_0-{\mathrm e}^s \in (0,1)$. Then,
$y_0 -G{\mathrm e}^s \in (0,G)$ and, using the definition of $f_0$
in such an interval, we have that
  $$
  \mbox{$
  f_0(y_0-G{\mathrm e}^s)=-G^5f_0\big(\frac{y_0-G{\mathrm e}^s}{G} -1\big) =
-G^5f_0(y_0-G{\mathrm e}^s).
 $}
 $$
   The same argument applies in $(0,G)$ and in all the
intervals obtained in every step of this construction. This completes the
proof. $\qed$

 \smallskip

As we mentioned above, by reflection of the solutions constructed
above, we obtain analogous results for solutions with the
left-hand interface. We note that, in both cases (the left- and
right-hand interfaces), the behavior of the self-similar profile
$F$ of (\ref{2.3}) close to the interfaces is described by means
of $f_0(y;G_2)$. The  behavior corresponds to TW ``travelling'' in
the opposite directions with $\l=+1$ is given by the profile
$f_0(y;G_1)$.

\subsection{Uniqueness of oscillatory TW solutions}

In Section \ref{sectIVP}, we have proved uniqueness of the
solution to the Cauchy problem for the equation (\ref{2.01}) with
given initial conditions (\ref{ivp1c}),
when $\alpha_j \neq 0$ at least for some values of $j$. It is
clear, due to the regularity conditions at the interface $y_0$,
that this analysis of uniqueness does not apply to the oscillatory
travelling wave solution at this point. In fact, it is not
difficult to check that uniqueness fails and, for instance, the
solutions $f \equiv 0$, the positive solution constructed above,
and its negative counterpart are three different solutions of the
same problem with analogous conditions at the interface $y=y_0$.

\ssk

In order to understand the structure of the set of solutions of
this type, it would be interesting to fix an additional condition
that provides uniqueness of the TW. We next prove that uniqueness
holds by adding to the previous conditions at the interface $y_0$,
a fixed transversal zero $a_0$ of the solution and a sign to the
function close to the fixed zero. Without loss of generality, we
next assume $a_0 =0$ and positivity of the solution to the
left-hand side of $a_0$. The result is stated as follows:

\begin{theorem}
 \label{Th.un1}
There exists a unique solution of the  equation $(\ref{b1})$ with
the interface $y_0$ with $f^{(j)}(y_0)=0$ for $j=0,1,2,3,4$ and
satisfying $f(0)=0$ and $f$ positive in $(-\e,0)$. Moreover, the
solution is given by $f_0(y)$, with the algebraic construction
above.

 \end{theorem}

 \noi{\em Proof.}
Let $f(y)$ be a solution of  (\ref{b1}) satisfying the stated
assumptions. Define the auxiliary  function $g(y) = f(y) + p(y)$,
with $p(y)$ a polynomial of the fourth degree,
 $$
  \mbox{$
 p(y)=\sum\limits_{i=1}^4 \frac{m_j y^j}{j!}\, ,
 $}
  $$
   with
$m_j=(f_0-f)^{(j)}(0)$. By the definition of $g$ and taking into
account the properties of $f$ and $f_0$ around $a_0=0$, it is not
difficult to check that $g(y)$ satisfies the same equation as
$f_0$ on $(-\e, \e)$ and that $g^{(j)}(0)=f^{(j)}(0)$ for
$j=0,1,...,4$.

 Hence, it follows from the uniqueness result of the solutions to
the initial value problem stated in Section \ref{sectIVP}, that $g
\equiv f_0$ at least in this interval, and, by a continuation
argument (applied for instance to the integral version of the
equation), we conclude the equivalence of both functions and its
derivatives for every $y \le y_0$.

\ssk

  In particular, we have at the interface $y_0$ that $g^{(j)}(y_0)=0$
for $j=0,1,2,3,4$. Taking into account the definition of $g$ and that
$f_0$ and $f$ satisfy also these regularity conditions at $y_0$, it
follows that $p^{(j)}(y_0) =0$ for $j=0,1,2,3,4$, whence $p(y) \equiv 0$ and
the uniqueness result follows. $\qed$

\ssk

This  allows us to establish a result concerning the structure of
the set of oscillatory TW solutions to equation (\ref{b1}) with
finite interfaces.

\begin{theorem}
 \label{Th.un2}
Assume that $f(y)$ is an oscillatory TW solution of $(\ref{b1})$ with  a finite
interface and necessary regularity assumptions  at the interface. Then, there
exist real values $\gamma >0 $ and $D$ such that
$$
 \mbox{$
f(y)= \pm \gamma^5 f_0(\frac {y}{\gamma} + D).
 $}
$$
\end{theorem}

 \noi{\em Proof.} By using the rescaling properties of the solutions,
 it is enough to define an auxiliary function
 $$
  \mbox{$
 g(y) = \pm \lambda^5 f(\frac y{\lambda} + C)
  $}
 $$
 with
 $\lambda >0$ and $C$ such that $y_0$ is the interface of $g(y)$ and
 $g(0)=0$. We also choose the sign in this definition in order to
 have all the hypotheses in Theorem \ref{Th.un1} to hold, whence we deduce that $g (y) \equiv f_0 (y)$. It is not difficult
to check that the result follows for $\gamma=\frac 1\lambda$ and
$D=- \frac C\lambda$. $\qed$

\ssk

Both results have strong implications concerning the properties
and the character of the oscillations of the solutions and
therefore, the behaviour of solutions close to the interface. In
particular, we  prove the following:

\begin{proposition}
 \label{limsup}
Assume that $f(y)$ is an oscillatory TW solution with the
interface $y_0$. Then there exists a finite limit:
   $$
 \mbox{$
\limsup\limits_{y \to y_0} \frac {f(y)}{(y_0-y)^5} =C.
 $}
$$
\end{proposition}
\noi{\em Proof.} We prove the result for $f_0(y)$. Denote by $y_1
\in (-1,0)$, the absolute maximum point of $f_0(y)$ in the
interval $(-1,0)$ and define
  $$
   \mbox{$
  y_{n+1} = (y_n + 1)G, \quad \mbox{with} \quad  G=\frac {y_0}{1 + y_0}.
   $}
  $$
   From the properties of
$f_0$ and its construction, it is not difficult to see that, so
defined, $y_n$ are the absolute maximum point of $|f_0(y)|$ in
every interval $I_n$ of the extension defined in the algebraic
construction of $f_0$. Taking into account the definition of $f_0$
and $y_n$ yields, as $n \to \iy$,
 $$
  \mbox{$
 y_n= \frac{G -
(1-y_1)G^{n+1}}{1-G}
 $} \quad \mbox{and} \quad
 \mbox{$
\frac {|f(y_n)|}{(y_0-y_n)^5} = \frac {G^{5(n+1)}f(y_1)}{(\frac G
{1-G} - y_n)^5} \to C,
 $}
$$
  so the result follows. $\qed$

\ssk

As a straightforward consequence of the construction above, we
also obtain the structure of the zeros of the solution $f_0$
constructed:

\begin{corollary}
 \label{Cor.Tr}
Every zero of the solution $f_0(y)$, except the one corresponding
to the interface, is transversal, i.e., $f' \not =0$.
\end{corollary}




\section{Continuity geometric extension to $m \approx 1$}
 \label{SMap1}

As the next natural step, we consider  equation (\ref{2.5r}) with
$\a \approx 0$. We then are going to perform the same geometric
(but not algebraic!) construction explained in the steps
(i)--(iii) in Section \ref{S6.1}, where the only difference is
that the basic profile $f_0(y)$ is not given explicitly as in
(\ref{b2}). We can use the standard  continuous dependence of
local solutions of the ODE (\ref{2.5r}) on the parameter $\a
\approx 0$. Here, we mean a class of good locally and uniquely
extensible solutions without strong singularities and finite
interface points.
 Of course, we crucially need the zero transversality property fixed in Corollary \ref{Cor.Tr}.
Good local properties of such solutions are checked by standard
contractivity techniques (as in Section \ref{S3.2}), so we arrive
at a local continuity (``homotopy $\a$-deformation") result:

\begin{theorem}
 \label{Th.Exm2}


 Equation $(\ref{2.71})$ with both signs $``\pm"$ $($i.e., for $\l= \mp 1)$ has a nontrivial periodic
 solution $\var_*(s)$ of changing sign for all  $\a \approx 0$.

 \end{theorem}

 It is worth recalling here that, in a global sense,  for $\a$ sufficiently far away from 0,
 periodic solutions do not exist; cf. Proposition \ref{Pr.Non} and
the heteroclinic bifurcation phenomena  in (\ref{mh02}).

Notice that the existence result in Theorem \ref{Th.Exm2} is not associated with the
``hyperbolicity" of the periodic orbit $\varphi_*$ in this range.
Anyway, Theorem \ref{Th.Exm2} altogether  with Propositions
\ref{Pr.Non} and \ref{Ptr.Stab}, though not covering the whole
$\a$-range, provide us with rather solid mathematical evidence on
existence of periodic solutions and also on the type of
heteroclinic bifurcations, at which these disappear.





\mpb

\noindent{\bf Acknowledgement. }  The first author thanks the
Department of Mathematical Sciences of the University of Bath for
its hospitality during her visits and Projects MTM2007-65018 (MEC)
and SA104A06 (JCyL) for financial support.



\bibliographystyle{amsplain}

\begin{thebibliography}{10}





\bibitem
 {Bern88}
 F.~Bernis, \emph{Source-type solutions of fourth order
 degenerate parabolic equations}, In: {Proc. Microprogram
 Nonlinear Diffusion Equation and Their Equilibrium States},
  W.-M.~Ni, L.A.~Peletier, and J.~Serrin, Eds., MSRI Publ., Berkeley,
 California, Vol. 1, New York, 1988, pp. 123--146.


  \bibitem
{Bern01}
  F.~Bernis, {\em Finite speed of propagation and
asymptotic rates for some nonlinear higher order parabolic
equations with absorption}, {Proc. Roy. Soc. Edinburgh}, {\bf
104A} (1986), 1--19.





\bibitem 
 {BF1}
 F.~Bernis and A.~Friedman, \emph{Higher order nonlinear degenerate
 parabolic equations}, J.~Differ. Equat., \textbf{83} (1990), 179--206.


\bibitem 
 {BMc91}
 F.~Bernis and J.B.~McLeod, \emph{Similarity solutions of a higher order
nonlinear
 diffusion
 equation}, Nonl. Anal., TMA, \textbf{17} (1991), 1039--1068.





\bibitem
 {CGII}
 M.~Chaves and V.A.~Galaktionov, {\em  On 
 source-type solutions  and the Cauchy problem
for a doubly degenerate sixth-order thin film
 equation II. Global properties}, in preparation
 (in arXiv.org shortly).


\bibitem{EidSys}
S.D.~Eidelman, {Parabolic {S}ystems}, North-Holland Publ. Comp.,
  Amsterdam/London, 1969.










   \bibitem{Gl4}
J.D.~Evans, V.A.~Galaktionov, and J.R.~King, \emph{Source-type
solutions of the fourth-order unstable thin film equation}, Euro
J.~Appl. Math., {\bf 18} (2007), 273--321.


      \bibitem
      {GBl6}   
J.D.~Evans, V.A.~Galaktionov, and J.R.~King, {\em Unstable
sixth-order thin film equation. I. Blow-up similarity solutions;
II. Global similarity patterns},
 {Nonlinearity}, {\bf 20} (2007), 1799--1841, 1843--1881.















\bibitem
 {GalRDE4n}
 V.A.~Galaktionov,
 {\em Countable branching of similarity solutions of higher-order
 porous medium type equations}, Adv. Differ. Equat., {\bf 13}
 (2008), 641--680.


 \bibitem
 {GalNDE5}
 V.A.~Galaktionov,
 {\em Shock waves and compactons for fifth-order nonlinear dispersion
 equations}, Europ. J.~Appl. Math., 2009, to appear (arXiv:0902.1632v1).



 \bibitem
  {PetI}
V.A.~Galaktionov and P.J.~Harwin, {\em On centre subspace
behaviour in thin film equations}, SIAM J.~Appl. Math., {\bf 69}
(2009), 1334--1358
 (an earlier preprint in arXiv:0901.3995v1).


\bibitem 
{GS1S-V} V.A.~Galaktionov and A.E.~Shishkov, {\em Saint-Venant's
principle in blow-up for higher-order quasilinear parabolic
 equations}, {Proc. Royal Soc. Edinburgh, Sect.~A}, {\bf 133} (2003), 1075--1119.




\bibitem
 {GSVR} V.A.~Galaktionov and S.R.~Svirshchevskii, Exact Solutions and
 Invariant Subspaces of Nonlinear Partial Differential Equations in Mechanics and Physics,
  Chapman$\,\&\,$Hall/CRC, Boca Raton,
Florida,
 2007. 





      \bibitem
      {Gia08}
L.~Giacomelli, H.~Kn\"upfer, and F.~Otto, {\em Smooth
zero-contact-angle solutions to a  thin film equation
 around the steady state},
 {J.~Differ. Equat.,} {\bf 245} (2008), 1454--1506.


\bibitem 
{Grun04} G.~Gr\"un, {\em Droplet spreading under weak slippage --
existence for the Cauchy problem}, {Commun. Partial Differ.
Equat.,} {\bf 29} (2004), 1697--1744.






\bibitem 
{Ki01} J.R.~King, {\em Two generalisations of the thin film equation},
{Math. Comput. Modelling,} {\bf 34} (2001), 737--756.




\bibitem{KrasZ}
M.A.~Krasnosel'skii and P.P.~Zabreiko, {Geometrical Methods of
Nonlinear
  Analysis}, Springer-Verlag, Berlin/Tokyo, 1984.



\bibitem 
{Perko} L.~Perko, {\rm Differential Equations and Dynamical
Systems}, Springer-Verlag, New York, 1991.



\bibitem 
 {RosH93}  P.~Rosenau and J.M.~Hyman, {\em Compactons: solitons with
 finite wavelength},
 {Phys. Rev. Lett.,} {\bf 70} (1993), 564--567.



\bibitem
{Shi2}
 A.E.~Shishkov, {\em Dead cores and instantaneous
compactification of the supports of energy solutions of
quasilinear parabolic
equations of arbitrary order}, {Sbornik: Math.}, {\bf 190}    
 (1999), 1843--1869.





\bibitem 
{VainbergTr} M.A.~Vainberg and V.A.~Trenogin, {\rm Theory of
Branching of Solutions of Non-Linear Equations}, Noordhoff Int.
Publ., Leiden, 1974.


  \bibitem 
  {Wu01}
   Z.~Wu, J.~Zhao, J.~Yin, and H.~Li, {Nonlinear Diffusion
  Equations,} World Scientific Publ. Co., Inc., River Edge, NJ, 2001.

\end{thebibliography}


\end{document}